%% file: main.tex
\documentclass[a4paper]{amsart}
\usepackage[a4paper]{geometry}
\usepackage[latin1]{inputenc}
\usepackage[T1]{fontenc}
\usepackage{lmodern}
\usepackage{amssymb,amsmath,amsthm,amscd}
\usepackage{mathrsfs}
\usepackage{graphicx} %
\usepackage{hyperref}
\usepackage{nameref,color,url,thmtools,thm-restate,zref-user,cleveref,stmaryrd,ifthen,mathtools,fancyvrb,bbm}
\usepackage{url,stmaryrd,ifthen}
\usepackage[all]{xy}
\usepackage{chngcntr}
\usepackage{marginnote}
\usepackage{diagbox,pict2e,tabu}
\usepackage{microtype}
\usepackage{mparhack}
\usepackage{etoolbox}
\usepackage[all]{xy}
\DeclareFontShape{OT1}{cmtt}{bx}{n}{<5><6><7><8><9><10><10.95><12><14.4><17.28><20.74><24.88>cmttb10}{}
\definecolor{darkblue}{RGB}{0,0,96}
\definecolor{gray}{RGB}{127,127,127}
\definecolor{darkgreen}{RGB}{0,160,0}
\definecolor{darkred}{RGB}{255,0,0}
\hypersetup{colorlinks={true},linkcolor={darkblue},citecolor={darkblue},filecolor={darkblue},urlcolor={darkblue},pdfauthor={Lukas Lewark},pdftitle={Rasmussen's spectral
sequences and the slN concordance invariants}}
\declaretheorem[name=Theorem,refname={theorem,theorems}]{Thm}
\declaretheorem[numberwithin=section,name=Proposition,refname={proposition,propositions}]{Prop}
\declaretheorem[sibling=Prop,name=Corollary,refname={corollary,corollaries}]{Cor}
\declaretheorem[sibling=Prop,name=Lemma,refname={lemma,lemmas}]{Lem}
\declaretheorem[sibling=Prop,name=Remark,refname={remark,remarks},style=remark]{Rmk}
\declaretheorem[sibling=Prop,name=Example,refname={example,examples}]{Ex}
\declaretheorem[numbered=no,name=Definition]{Def}
\declaretheorem[numbered=no,name=Conjecture]{Conj}
\newenvironment{Pf}[1][Proof]{\begin{proof}[#1]}{\end{proof}}
\newcommand{\thmref}[1]{\hyperref[#1]{\namecref{#1}~\ref{#1}}}
\newcommand{\Thmref}[1]{\hyperref[#1]{\nameCref{#1}~\ref{#1}}}
\newcommand{\dover}[1]{\overline{\overline{#1}}}

\DeclareMathOperator{\tr}{tr}
\DeclareMathOperator{\xdeg}{xdeg}
\newenvironment{renumerate}{\begin{enumerate}}{\end{enumerate}}

\newcommand{\cat}[1]{\ensuremath{\mathsf{#1}}}
\newcommand{\Red}[1]{\overline{{#1}}}

\newcommand{\emptypoint}{\ensuremath{\,\cdot\,}}

\newlength{\lowerhalftmp}
\newcommand{\aFigpath}{figs}
\newcommand{\lowerhalfx}[1]{\settoheight{\lowerhalftmp}{#1}\addtolength{\lowerhalftmp}{-1.2ex}\raisebox{-.5\lowerhalftmp}{#1}}

\newcommand{\avcfig}[2][]{\lowerhalfx{\includegraphics[#1]{\aFigpath/#2.pdf}}}
\newcommand{\afig}[2][]{{\includegraphics[#1]{\aFigpath/#2.pdf}}}
\newcommand{\homf}{\textsc{Homflypt}}
\newcommand{\Sln}{\ensuremath{\mathfrak{sl}_N}}
\newcommand{\CC}[2]{\ensuremath{C_{#2}(#1)}}
\DeclareMathOperator{\xdim}{xdim}
\newcommand{\Nat}{\ensuremath{\mathbb{N}}}
\newcommand{\Int}{\ensuremath{\mathbb{Z}}}

\newcommand{\Hom}[3][]{\ifthenelse{\equal{#1}{rr}}{\ensuremath{\dover{\left\llbracket #2\right\rrbracket}_{#3}}}{\ifthenelse{\equal{#1}{r}}{\ensuremath{\overline{\left\llbracket #2\right\rrbracket}_{#3}}}{\ensuremath{\left\llbracket #2\right\rrbracket_{#3}}}}}

\newcommand{\qua}{\hskip 0.75em plus 0.15em \ignorespaces}
\def\arxiv#1{\relax\ifhmode\unskip\qua\fi
    \href{http://arxiv.org/abs/#1}%
{\tt arXiv:\penalty -100\unskip#1}}    

\def\MR#1{\relax\ifhmode\unskip\qua\fi
    \href{http://www.ams.org/mathscinet-getitem?mr=#1}{MR#1}}
\def\xox#1{\csname xx#1\endcsname}

\begin{document}
\title{Rasmussen's spectral sequences and the \Sln-concordance invariants}
\author{Lukas Lewark}
\address{Mathematical Sciences\\
Durham University\\
UK}
\email{lukas.lewark@durham.ac.uk}
\urladdr{http://www.maths.dur.ac.uk/$\sim$vxhn54}
\date{November 21, 2013}
\thanks{Supported by the EPSRC-grant EP/K00591X/1.}
\begin{abstract}
\input{abstract}
\end{abstract}
\maketitle
\setcounter{tocdepth}{1}
\setcounter{secnumdepth}{1}
\tableofcontents
\bibliographystyle{myamsalpha}
\section{Introduction}
\input{intro}
\bigskip\paragraph{\emph{Acknowledgements}}
This paper is extracted from my thesis \cite{these},
and I thank Christian Blanchet for having been my adviser.
Thanks to Andrew Lobb for comments on a first version of the paper.
\section{The Khovanov-Rozansky homologies}
\label{sec:overview}
\input{overview}
\section{The reduced-unreduced spectral sequence}
\label{sec:sseq}
\input{sseq}
\section{From \homf-homology to the \Sln-concordance invariants}
\label{sec:tool}
\input{tool}
\section{Slice-torus knot concordance invariants}
\label{sec:slicetorus}

\input{slicetorus}
\section{Linear independence of some of the \Sln-concordance invariants}
\label{sec:examples}
\input{examples}
\newpage
\bibliography{main}
\end{document}

%% file: abstract.tex
Combining known spectral sequences
with a new spectral sequence relating reduced and unreduced \Sln-homology
yields a relationship between the \homf-homology of a knot and its \Sln-concordance invariants.
As an application, some of the \Sln-concordance invariants are shown to be linearly independent.

%% file: intro.tex
The Khovanov-Rozansky homologies \cite{roz,roz2} are categorifications of the \homf-polynomial and
its specialisations.
There is a wealth of different homology theories.
As an astonishing and powerful application, some of them induce concordance invariants which give lower bounds
for the smooth slice genus of a knot \cite{rasmussen,wu4,lobbonroz,lobbconcordance}.
This paper has two goals: to explore how the different Khovanov-Rozansky
homologies are related; and to show that the said concordance invariants, though very
close to each other, are not equal.
Apart from trying to deepen the understanding of the Khovanov-Rozansky homologies themselves,
there is a geometric motivation: for example, the Khovanov-Rozansky concordance invariants
may detect free summands in the group of topologically slice modulo smoothly slice knots \cite{slicelivingston};
and they could even be used to disprove the smooth Poincaré conjecture \cite{manandmachine}.
The main result is the following:
\begin{Thm}\label{c1}
Let $\tau$ be the concordance invariant from knot Floer homology \cite{os,rasphd}, and for all $N \geq 2$,
$s_N$ the concordance invariant from $\Sln$-homology (see \thmref{prop:sndef}).

\begin{renumerate}
\item Neither $\tau$ nor $s_2$ is a linear combination of $\{s_N\}_{N\geq 3}$, and for any fixed $N\geq 3,
\{\tau, s_2, s_N\}$ are linearly independent. (Proof without computer calculations.)
\item %
$s_3$ is not a linear combination of $\{s_N\}_{N\geq 4}$, and for any fixed $N\geq 4,
\{\tau, s_2, s_3, s_N\}$ are linearly independent. (Proof relies on computer calculations.)
\end{renumerate}
\end{Thm}
This prompts the following conjecture:
\begin{Conj}
The concordance invariants $\{\tau\}\cup \{s_i\}_{i \geq 2}$ are
linearly independent.
\end{Conj}

See \cite{hom,jabuka} for similar results in Heegard-Floer homology.

The Khovanov-Rozansky homologies considered in this paper are:
the triply graded homology $\Hom[r]{\emptypoint}{\infty}$
 categorifying the \homf-polynomial;
for each $N \geq 1$, the doubly graded reduced homology
 $\Hom[r]{\emptypoint}{N}$
 and unreduced homology 
 $\Hom{\emptypoint}{N}$
categorifying
the reduced and unreduced \Sln-polynomial, respectively; and the deformation of the unreduced
homology, the filtered homology
 $\Hom{\emptypoint}{N}^f$.
On the uncategorified level, the different polynomials are specialisations or multiples with a fixed factor of one another.
On the categorified level, relations normally take the form of spectral
sequences; but the understanding of the interdependence of the different Khovanov-Rozansky
homologies is far from complete. The following theorem clarifies the relationship between unreduced
and reduced \Sln-homology.
Let $D$ be a diagram of a link $L$ with a marked component, and $N \geq 1$ an integer.
Let $\CC{D}{N}$ be the graded chain complex
defined by Khovanov and Rozansky \cite{roz} whose homology is $\Hom{L}{N}$.
\begin{Thm}\label{thm:sseq}
There is a filtration of $C_N(D)$ respected by the differential such that
the induced spectral sequences satisfy the following properties
(where $r$ denotes the grading associated to the new filtration):
\begin{renumerate}
\item Its differentials respect the $q$-degree; it converges on the $N$-th page, and forgetting the $r$-grading,
the limit is isomorphic to $\Hom{L}{N}$.
\item Its first page is isomorphic to
\[
\frac{(qr)^N - (qr)^{-N}}{qr - (qr)^{-1}}\cdot \Hom[r]{L}{N}.
\]
\item The higher pages are invariants of links with a marked component.
\end{renumerate}
\end{Thm}
Let us have a closer look at the Khovanov-Rozansky concordance invariants.
They belong to the broader class of slice-torus invariants introduced (without a name) by Livingston \cite{livingston}:
\begin{Def}\label{def:slice torus}
A \emph{slice-torus knot invariant} is a homomorphism $\nu$ from the smooth knot concordance group
to the real numbers that satisfies the following conditions:
\begin{itemize}
\item[\textbf{(slice)}]
For all knots $K$, $\nu(K)$ is a lower bound to twice the slice genus: $\nu(K) \leq 2g_4(K)$.
\item[\textbf{(torus)}]
For positive torus knots, this bound is sharp, i.e.
\[
\forall p, q\in\mathbb{Z}^+: (p,q) = 1 \implies\nu(T(p,q)) = 2g_4(T(p,q)) = (p-1)(q-1).
\]
\end{itemize}
\end{Def}
Note that we chose a different normalisation than Livingston. More importantly, we consider
real-valued instead of integer-valued invariants, in order to include the normalised
Khovanov-Rozansky concordance invariants. However, every slice-torus invariant discussed in this text
takes values only in $\frac{1}{n}\mathbb{Z}$ for some fixed $n$.

Slice-torus knot invariants form a closed convex subset of the space of all
real concordance homomorphisms.
The slice-torus conditions are quite restrictive;
all slice-torus invariants have e.g. the same value on quasi-positive knots.
In this paper, the sharper slice-Bennequin inequality is generalised to slice-torus invariants,
thereby showing that all slice-torus invariants have the same value on homogeneous knots
(see \thmref{lem:lobbineq} in \thmref{sec:slicetorus}).

Up to convex linear combination, $2\tau$ (shown to be distinct from the Rasmussen invariant in \cite{sneqtau})
is the only known slice-torus invariant
not stemming from the Khovanov-Rozansky homologies.
The oldest Khovanov-Rozansky concordance invariant
is the Rasmussen invariant $s_2$ \cite{rasmussen}, which comes from
Khovanov homology over a field of characteristic 0.
Later on, generalisations were defined: an invariant $s_N$ coming from 
$\mathfrak{sl}_N$-homology for arbitrary $N \geq 2$ \cite{lobbconcordance};
and an invariant $s_2^{\mathbb{F}_p}$, obtained from $\mathfrak{sl}_2$-homology
over a prime field $\mathbb{F}_p$
\cite{tanglescobordisms, turnermod2} (see also \cite{remarkOnRas}).
So far, computer calculations indicate that $s_2 \neq s_2^{\mathbb{F}_2}$ \cite{lipsar,knotkit};
and that $s_3 \neq s_2$ \cite{lew1,foamho}.
The following theorem, whose proof uses \thmref{thm:sseq}, the Rasmussen spectral sequences \cite{somedifferentials}
and the Lee-Gornik spectral sequences \cite{lee,gornik}, is a means to distinguish the $s_N$.
%\pagebreak
\begin{Thm}\label{cor:tool}
Let $K$ be a knot. For all $N \geq 2$, let
\[
X_N = \{ \alpha + N\beta \mid 
\Hom[r]{K}{\infty} \text{ has a generator of degree } q^{\alpha} a^{\beta}\text{ in homological degree } 0 \} \subset 2\Int.
\]
Then for the unnormalised \Sln-concordance invariant $s_N'$, we have $\min X_N \leq s_N'(K)\leq \max X_N$, or equivalently,
for the normalised invariant:
\begin{equation*}
\frac{\max X_N}{1 - N} \ \leq\ s_N(K)\
 \leq\ 
\frac{\min X_N}{1 - N}. 
\end{equation*}
\end{Thm}

\Thmref{cor:tool} combined with the sharper slice-Bennequin inequality enables
us to calculate the \Sln-concordance invariants of certain three-stranded pretzel knots:
\begin{figure}[ht]
\hfill\afig[scale=0.8]{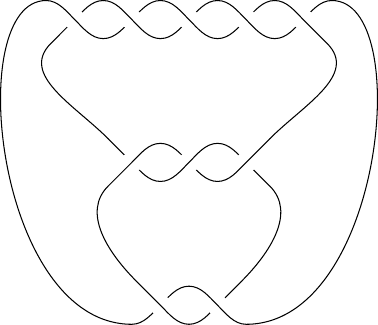}\hfill\ 
\caption[The $(5,-3,2)$-pretzel knot.]%
{The $(5,-3,2)$-pretzel knot $K$ with $s_2(K) = 2$, $s_3(K) = 1$ and $s_N(K) \in \{0, 2/(N-1)\}$ for all $N \geq 4$.}
\label{fig:pretzel}
\end{figure}
\begin{Thm}\label{cor:examples}
Let $\ell$ and $m$ be odd integers and $\ell > m \geq 3$. Then
\begin{align}
\tag{\textbf{i}}
s_2(P(\ell, -m, 2)) & = \ell - m,\quad\text{and}\\
\tag{\textbf{ii}}
\forall N \geq 3:
s_N(P(\ell, -m, 2)) & \in \Bigl\{\ell - m - 2, \ell - m - 2 + \frac{2}{(N-1)}\Bigr\}.
\end{align}
Let $n \geq 4$ be an even integer. Then
\begin{align}
\tag{\textbf{iii}}
s_2(P(\ell, -m, n)) & =
\begin{cases}
\ell - m & \text{if } m > n\\
\ell - m - 2 & \text{if } m < n,\qquad\text{and}
\end{cases} \\
\tag{\textbf{iv}}
\forall N \geq 3:
s_N(P(\ell, -m, n)) & = \ell - m - 2.
\end{align}
\end{Thm}
These pretzel knots are the examples used to prove \thmref{c1}.
The remainder of the paper is organised as follows:
\thmref{sec:overview} details the needed results on the Khovanov-Rozansky homologies.
So while this paper is de facto self-contained, nevertheless some familiarity with
\cite{roz} is advisable.
The two following sections contain the proofs of
\thmref{thm:sseq} and \thmref{cor:tool}, respectively.
In \thmref{sec:slicetorus}, known results about slice-torus invariants are collected,
and the sharper slice Bennequin inequality (\thmref{lem:lobbineq}) is proven.
\Thmref{sec:examples} finally applies the tools to
examples such as the family of pretzel knots and contains the proofs of
\thmref{cor:examples} and \thmref{c1}.

%% file: overview.tex
This section gives an overview over the different Khovanov-Rozansky homologies
and the known spectral sequences relating them. Notations and conventions are clarified on the way.
Comparing Khovanov-Rozansky homologies over different base fields
leads to interesting open questions (cf. \cite{lipsar,remarkOnRas,turnermod2}).
However, we restrict our attention to characteristic 0, and
consider all chain complexes to be over the complex numbers.

\subsection{Unreduced homology}
Let $D$ be a diagram of a link $L$.
For all integers $N \geq 1$, Khovanov and Rozansky \cite{roz}
define a chain complex $\CC{D}{N}$ of graded vector spaces (technically, it is a cochain complex).
Any Reidemeister move gives rise to a quasi-isomorphism of
$\CC{D}{N}$, and so the homology $\Hom{L}{N}$ is a link invariant,
called the \Sln-homology. Note that all links have the same $\mathfrak{sl}_1$-homology,
while $\mathfrak{sl}_2$-homology is isomorphic to Khovanov homology \cite{khovanov},
and $\mathfrak{sl}_3$-homology to a homology theory defined previously via webs and foams \cite{khovanovsl3,vaz}.

We regard this chain complex and his homology as a doubly graded vector space, with
a \emph{homological} ($t$), and a \emph{quantum} ($q$) degree. In general, for such a graded space $V$,
let us write $V^i$ for the subspace of homological degree $i$,
and, if $V$ has finite dimension, $\xdim V\in \Nat[t^{\pm 1}, q^{\pm 1}]$ (where $\Nat = \{0,1,2,\ldots\}$) for the graded dimension of $V$.
If $f$ is a homogeneous homomorphism of such spaces and has
$(t,q)$-degree $(i,j)$, we denote by $\xdeg f$ the monomial $t^i q^j$.

The \Sln-homology categorifies
the \Sln-polynomial $P_N$ \cite{RT}, i.e.
$\xdim(\Hom{L}{N})(-1,q) = P_N(L)$. The \Sln-polynomial 
is given by
its value of $[N]_q = (q^{-N+1} + q^{-N+3} + \ldots + q^{N-1})$
on the unknot
and
the following skein relation: %
\begin{equation*}
q^N\cdot P_N\Bigl(\avcfig{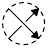}\Bigr) -
q^{-N}\cdot P_N\Bigl(\avcfig{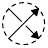}\Bigr) =
(q - q^{-1})\cdot P_N\Bigl(\avcfig{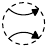}\Bigr).
\end{equation*}

In fact, the \Sln-homology theory is richer than that, being defined for tangles as well:
to every tangle diagram $D$ with boundary $\partial D$ (the boundary being a finite sequence
of signs) an object in some category $\cat{C}_{\partial D}$
is associated.
This is done in such a way that the gluing of tangle diagrams corresponds to
the tensor product of the associated objects (i.e. this is a \emph{canopolis}-construction, cf. \cite{tanglescobordisms}).
The category $\cat{C}_{\varnothing}$ is 
equivalent to the category of graded chain complexes over $\mathbb{C}$.
In this way, Reidemeister invariance of the \Sln-homology of a link can be proven
simply by showing that the objects associated to the two small tangle diagrams
which correspond to each of the Reidemeister moves are isomorphic.

\subsection{Reduced homology}
The reduced version of this homology
categorifies the reduced \Sln-polynomials $\Red{P_N} = P_N / [N]_q$.
Let $D$ be a diagram of the link $L$ with a marked component. Let $A = \mathbb{C}[X]/(X^N)$,
an algebra with grading $\deg X^i = 2i$ for $i\in\{0,\ldots, N - 1\}$.
Then $\CC{D}{N}$ has the structure of a free graded $A$-module.
This structure is respected by the differential of $\CC{D}{N}$, and it may depend on the choice of the marked component.
Let $\widetilde{\mathbb{C}}$ be the graded $A$-module $A / (X)$ with a shift of $N-1$ in the $q$-grading.
Let $\Red{\CC{D}{N}} = \CC{D}{N} \otimes_A \widetilde{\mathbb{C}}$.
The following proposition, which is essential for the proof of \thmref{thm:sseq}(iii),
 is implicit in \cite[end of section 7]{roz}. Let us give an explicit proof.
\begin{Prop}\label{prop:red}
If two base-pointed diagrams $D$ and $D'$ are connected by a Reidemeister move that avoids
the base-point, then there is a chain homotopy equivalence respecting the $A$-module structure
between $\CC{D}{N}$ and $\CC{D'}{N}$.
\end{Prop}
\begin{Pf}
This proof supposes greater familiarity with the details of the construction \cite{roz} than the rest of the paper:
in fact, we consider link diagrams $D$ with \emph{marks}.
Marks form a finite subset of $D$ that avoids the crossings, such that any interval connecting
two crossings contains at least one mark. The complement of the marks is the disjoint union of
components, all of which are either a positive or negative crossing, a line or a circle.
The chain complex $\CC{D}{N}$ is then defined as the tensor product (over adequate rings)
of the elementary chain complexes associated to these simple pieces.
Adding or removing marks produces a homotopy equivalent chain complex, and this chain homotopy equivalence
respects the $A$-module structure (this follows from \cite[Proposition 22]{roz}).
So, without loss of generality, assume that $D = D_1 \cup D_3$ and $D' = D_2 \cup D_3$
each split along marks into two tangle diagrams:
small ones $D_1$ and $D_2$, in which the Reidemeister move takes place, who have the same complement $D_3$.
Then, $\CC{D}{N} = \CC{D_1}{N} \otimes \CC{D_3}{N}$ and $\CC{D'}{N} = \CC{D_2}{N} \otimes \CC{D_3}{N}$.
These tensor products are $A$-modules because the second factor is.
There is a chain homotopy equivalence $\varphi$ between $\CC{D_1}{N}$ and $\CC{D_2}{N}$.
The tensor product of $\varphi$ with the identity of $\CC{D_3}{N}$ gives a chain homotopy equivalence
between $\CC{D}{N}$ and $\CC{D'}{N}$ that respects the $A$-module structure.
\end{Pf}
\begin{Cor}[{\cite{roz}}]
The homology $\Hom[r]{L}{N}$ of the reduced complex is an invariant of links with a marked component.
\end{Cor}
\begin{Pf}
Two base-pointed diagrams $D$ and $D'$ represent the same base-pointed link if and only if they
are connected by a finite sequence of Reidemeister moves which avoid the base point.
\end{Pf}
\subsection{Filtered homology}
There is a filtered version of \Sln-homology, whose associated graded is the
original unreduced \Sln-homology.
As usual, a filtered complex gives rise to a spectral sequence.
Let us briefly clarify the indexing convention:
a spectral sequence is a sequence $(E_k)_{k\in\{0,1,2,\ldots\}}$ of graded chain complexes,
such that for all $k\geq 0$, forgetting the differential on $E_{k+1}$ yields
the homology of $E_k$.
If not specified otherwise, the differential $d_k$ on the $k$-th page has $(t,q)$-degree $(1, k)$;
this is non-standard, but convenient in our context.
\begin{Prop}[{\cite{lee} for $N = 2$, \cite{gornik} for all $N \geq 2$}]
\label{prop:Lees-spectral-sequence}
There is a spectral sequence starting at unreduced \Sln-homology
and converging to filtered \Sln-homology.
\end{Prop}
\begin{Prop}[{\cite{lee} for $N = 2$, \cite[Theorem 1.2]{wu4} for all $N \geq 2$}]
The higher pages of the Lee-Gornik spectral sequence are link invariants.
\end{Prop}
The following detail has not been explicitly stated:
\begin{Prop}
The differential on the $k$-th page of the
Lee-Gornik spectral sequence vanishes unless $k$ is a multiple of $2N$.
\end{Prop}
\begin{Pf}
Note that the differentials of Gornik's filtered chain complex
preserve the $q$-degree mod $2N$ (see \cite{gornik}).
So the chain complex decomposes as a direct sum of $N$ terms
(for $i\in\{0,\ldots, N-1\}$, the $i$-th term containing the generators
of $q$-degree equal to $2i$ mod $2N$),
and so does the induced spectral sequence.
\end{Pf}
As a consequence, it makes sense to forget all pages with vanishing differentials
and renumber: from now on, by the ``$k$-th page'' of the Lee-Gornik spectral sequences,
we actually refer to the $(2Nk)$-th page.
It is still an open conjecture that this spectral sequence converges on the second page.

\begin{Prop}[{\cite{rasmussen} for $N = 2$, \cite{gornik,wu4,lobbonroz,lobbconcordance} for all $N\geq 3$}]
\label{prop:sndef}
Let $K$ be a knot.
The filtered \Sln-homology of $K$ is isomorphic to the unreduced \Sln-homology of the
unknot, with a $q$-shift by some even integer which we denote by $s'_N(K)$.
Its normalisation $s_N(K) = s'_N(K) / (1 - N) \in \frac{2}{N-1}\mathbb{Z}$ is a slice-torus invariant called the
\Sln-concordance invariant.
\end{Prop}
Note that unlike the Rasmussen and twice the $\tau$-invariant, the \Sln-concordance invariants of a knot
need not be even integers, or integers at all.
\subsection{\homf-homology}
Khovanov and Rozansky \cite{roz2} introduce a chain complex $\CC{D}{\infty}$ of
doubly graded complex vector spaces defined for
a braid diagram $D$. Its homology %
is a link invariant called the \homf-homology, which categorifies the \homf-polynomial $P_{\infty} \in \mathbb{Z}[q^{\pm 1},a^{\pm 1}]$.
The \homf-polynomial is determined by
its value of $1$ on the unknot, and the following skein relation:
\[
\label{tag:homf skein relation}
a\cdot P_{\infty}\Bigl(\avcfig{35-pos}\Bigr) -
a^{-1}\cdot P_{\infty}\Bigl(\avcfig{32-neg}\Bigr)  =
(q - q^{-1})\cdot P_{\infty}\Bigl(\avcfig{03-Eqsign}\Bigr)
\]
There are several versions of \homf-homology. Rasmussen \cite{somedifferentials}
e.g. works with a reduced and an unreduced version, and an interpolation of the
two; but all these versions carry the same information (as is not the
case for the reduced and unreduced version of \Sln-homology).
In this text, we stick to the reduced version, denoted by $\Hom[r]{\,\cdot\,}{\infty}$.
For a knot this is, unlike the unreduced version, a finite dimensional space.
We follow similar grading conventions
as Mackaay and Vaz \cite{terasaka}, but exchanging $t$ and $t^{-1}$, i.e.
\[
\xdim\Hom[r]{{T(3,2)}}{\infty}
= a^{-2}q^2 + t^2a^{-2}q^{-2} + t^3a^{-4}.
\]
In \cite{somedifferentials}, Rasmussen follows still another grading convention; the monomial $q^ia^jt^k$
in that convention corresponds to the monomial $q^ia^jt^{(k-j)/2}$ in ours.

There is yet another version of \homf- and \Sln-homology, only defined for two-component links:
\emph{totally reduced homology}, denoted by $\Hom[rr]{\,\cdot\,}{\infty}$ and $\Hom[rr]{\,\cdot\,}{N}$, respectively.
We will not give its definition, because we only ever use the totally reduced homology of the positive Hopf link
(as calculated e.g. in \cite{terasaka}):
\[
\xdim \Hom[rr]{{T(2,2)}}{\infty} =
(a^{-1}q^{2} + ta^{-1} + t^{2}a^{-1}q^{-2}  + t^{3}a^{-3})\cdot t^{-1/2}.
\]
\subsection{The relationship of \homf-homology and reduced \Sln-homology}
The Rasmussen spectral sequences show that in a certain sense, \homf-homology
is the stabilisation of the \Sln-homologies as $N \to \infty$.
\begin{Prop}[{\cite{somedifferentials}}]\label{prop:Rasmussens-spectral-sequence}
Let $L$ be a link with a marked component.
For every $N \geq 1$, there is a spectral sequence with first page $\Hom[r]{L}{\infty}$.
Its limit is, after a regrading, isomorphic to the reduced \Sln-homology of $L$.
Explicitly, the regrading of the $(t,q,a)$-degree is $(i,j,\ell) \mapsto (i, j + N\ell)$.
The differential on the $k$-th page of the spectral sequence
has degree $tq^{2Nk}a^{-2k}$.
The higher pages are invariants of links with a marked component.
If $L$ is a knot, then for sufficiently large $N$, this sequence converges on the first page.
\end{Prop}
\subsection{Calculating \homf-homology}
It is in fact easier to calculate the \homf-homology of a knot than its \Sln-homology for some $N$.
See \cite{terasaka} for an exemplary calculation. Let us present the part of the tool-kit
which is necessary for the calculations in this paper.
\homf-homology is well-behaved under taking the connected sum:
\begin{Prop}[{\cite[Lemma 7.8]{somedifferentials}}]\label{lem:sum}
Let $L_1$ and $L_2$ be links, and $L_3$ any connected sum of $L_1$ and $L_2$. Then
$\displaystyle
\Hom[r]{L_3}{\infty} \cong \Hom[r]{L_1}{\infty} \otimes \Hom[r]{L_2}{\infty}.
$
\end{Prop}
\begin{Def}[\cite{rasbridge}]
Let the $\delta$-grading on $\Hom[r]{\emptypoint}{\infty}$ be defined by $\delta(t^iq^ja^k) = 2i + j + 2k$.
A knot is \emph{KR-thin} if its \homf-homology is supported in a single $\delta$-degree
that is equal to minus its signature.
\end{Def}
\begin{Lem}
The \homf-homology of a KR-thin knot $K$ is determined by its \homf-polynomial $P_{\infty}(K)$
and its signature $\sigma(K)$:
\[
\xdim \Hom[r]{L}{\infty} = (-t)^{-\sigma(K)/2}\cdot P_{\infty}(qt^{-1/2}, at^{-1}).
\]
\end{Lem}
\begin{Pf}
This immediately follows from the fact that $\xdim\Hom[r]{L}{\infty}(-1,q,a) = P_{\infty}(q,a)$.
\end{Pf}
\begin{Prop}[{\protect\cite{rasbridge},\protect\cite[Corollary 1]{somedifferentials}}]
Two-bridge knots are KR-thin.\label{prop:twothin}
\end{Prop}
\begin{Rmk}\label{rmk:quasialt}
Quasi-alternating links are a generalisation of alternating links introduced in \cite{branchedheegard}.
Quasi-alternating links have thin Khovanov and knot Floer homology \cite{quasialt},
and in particular the Rasmussen and twice the $\tau$-invariant of quasi-alternating knots equal their signature.
This can be proven via an unoriented skein relation. For $N\geq 3$, however, \Sln-homology
does not satisfy such a relation, 
and indeed there are even alternating knots which are not KR-thin \cite{somedifferentials}.
Still, the \Sln-concordance invariants of alternating knots (but not, in general, of quasi-alternating knots)
equal their signature for all $N$, since this is true of all slice-torus invariants (see \thmref{cor:alt}).
\end{Rmk}
\begin{figure}[h]
\hfill $K_+ = $ \avcfig{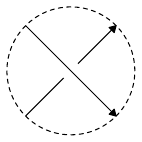}
\hfill $K_- = $ \avcfig{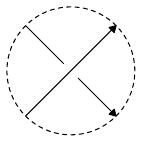}
\hfill $L_0 = $ \avcfig{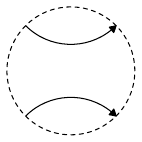}
\hfill\ 
\caption{The knots $K_{\pm}$ and the link $L_0$.}
\label{fig:Skein-Homf2}
\end{figure}
\begin{Prop}[The skein long exact sequences, {\cite[Lemma 7.6]{somedifferentials}}]
Let $K_+, K_-$ and $L_0$ be two knots and one two-component link which
look the same everywhere except near one crossing, where they differ as
shown in \thmref{fig:Skein-Homf2}. Then for all $N \geq 2$, there is a long exact sequence
\begin{equation*}
\cdots
\xrightarrow{\mathmakebox[1.3cm]{ }}
\Hom[r]{K_-}{N} \xrightarrow{\mathmakebox[1.3cm]{(-N,\frac{1}{2})}}
\Hom[rr]{L_0}{N} \xrightarrow{\mathmakebox[1.3cm]{(-N,\frac{1}{2})}}
\Hom[r]{K_+}{N} \xrightarrow{\mathmakebox[1.3cm]{(2N,-2)}}
\Hom[r]{K_-}{N} 
\xrightarrow{\mathmakebox[1.3cm]{ }}
\cdots 
\end{equation*}
The maps' $(t,q)$-degree is indicated above the arrows.
\end{Prop}
This proposition refers to \Sln-homology; to make a statement about \homf-homology, we need
the following technical lemma. Let $\leq$ denote the partial order of polynomials given
as follows: $A \leq B$ if and only if there is a polynomial $C$ with non-negative coefficients such that $A + C = B$.
\begin{Lem}
\label{lem:large N monomials}
Let $A, B \in \Nat[q^{\pm 1},a^{\pm 1}]$. Suppose that for infinitely many $N$,
$A(q,q^N) \leq B(q,q^N)$.  Then $A(q,a) \leq B(q,a)$.
\end{Lem}
\begin{Pf}
Let $i_{\text{max}}$ and $i_{\text{min}}$ be the maximal and minimal exponent of
$q$ occurring in $A$ or $B$. Choose $N$ such that $A(q, q^N) \leq B(q,q^N)$ and
$|N| > i_{\text{max}} - i_{\text{min}}$. Then different monomials in $A(q,a)$ and $B(q,a)$
yield different monomials in $A(q,q^N)$ and $B(q,q^N)$. To show this, consider
two monomials $c\cdot q^ia^j$ and $c'\cdot q^{i'}a^{j'}$ in $A(q,a)$ (with $c, c' \neq 0$).
Then $c q^{i + Nj} = c' q^{i' + Nj'}$ implies $c = c'$ and
$i + Nj = i' + Nj' \implies i - i' = N(j' - j) \implies |N|\cdot |j' - j| \leq i_{\text{max}} - i_{\text{min}}
\implies j' = j \implies i = i'$.
Let $c_{ij}$ and $c'_{ij}$ be the coefficients of the monomial $q^ia^j$
in $A(q,a)$ and $B(q,a)$, respectively. Then $c_{ij}$ and $c'_{ij}$ are also the respective coefficients of the monomial $q^{i + Nj}$
in $A(q,q^N)$ and $B(q,q^N)$, and thus $c_{ij} \leq c'_{ij}$.
\end{Pf}
\begin{Cor}
\label{prop:homf les}
Suppose $K_{\pm}$ and $L_0$ are given as in \thmref{fig:Skein-Homf2}, then
\begin{align*}
\xdim \Hom[r]{K_+}{\infty} & \quad\leq\quad  t^{2}\cdot a^{-2}\cdot \xdim \Hom[r]{K_-}{\infty} + t^{1/2}\cdot a^{-1} \cdot \xdim \Hom[rr]{L_0}{\infty}\quad\text{and} \\
\xdim \Hom[r]{K_-}{\infty} & \quad\leq\quad  t^{-2}\cdot a^{2}\cdot \xdim \Hom[r]{K_+}{\infty} + t^{-1/2}\cdot a \cdot \xdim \Hom[rr]{L_0}{\infty}.
\end{align*}
\end{Cor}
\begin{Pf}
We will just prove the first equation, the second one follows similarly.
The long exact sequence can be broken up into short ones; i.e., for some
quotient space $A$ of $\Hom[rr]{L_0}{N}$ and subspace $B$ of $\Hom[r]{K_-}{N}$
there is a short exact sequence
\[
0
\xrightarrow{\mathmakebox[1.3cm]{ }}
A
\xrightarrow{\mathmakebox[1.3cm]{(-N,\frac{1}{2})}}
\Hom[r]{K_+}{N} \xrightarrow{\mathmakebox[1.3cm]{(-2N,2)}}
B
\xrightarrow{\mathmakebox[1.3cm]{ }}
0.
\]
This is equivalent to $\Hom[r]{K_+}{N} \cong (q^{-N}t^{1/2}\cdot A) \oplus (q^{-2N}t^2\cdot B)$.
In terms of graded dimensions, this implies
\begin{align*}
\xdim\Hom[r]{K_+}{N} & = q^N\cdot t^{-1/2} \cdot \xdim A + q^{2N}\cdot t^{-2}\cdot\xdim B \\
\implies
\xdim\Hom[r]{K_+}{N} & \leq q^N\cdot t^{-1/2} \cdot \xdim \Hom[rr]{L_0}{N} + q^{2N}\cdot t^{-2}\cdot\xdim \Hom[r]{K_-}{N}.
\end{align*}
For large enough $N$, the three polynomials in this inequality stabilise (see \thmref{prop:Rasmussens-spectral-sequence}), i.e.
\begin{align*}
(\xdim\Hom[r]{K_+}{\infty})(q, q^N) &  = (\xdim\Hom[r]{K_+}{N})(q), \\
(\xdim\Hom[rr]{L_0}{\infty})(q, q^N) & = (\xdim\Hom[rr]{L_0}{N})(q), \\
(\xdim\Hom[r]{K_-}{\infty})(q, q^N) &  = (\xdim\Hom[r]{K_-}{N})(q).
\end{align*}
So using \thmref{lem:large N monomials}, the statement follows.
\end{Pf}
\begin{Rmk}
It may seem cumbersome to prove \thmref{prop:homf les} by first considering the \Sln-homologies, and then
\homf-homology as their stabilisation. But as Rasmussen remarks, who proves the KR-thinness of two-bridge knots
in the same way \cite{somedifferentials}, it is unclear how to work directly on \homf-homology.
\end{Rmk}

%% file: sseq.tex
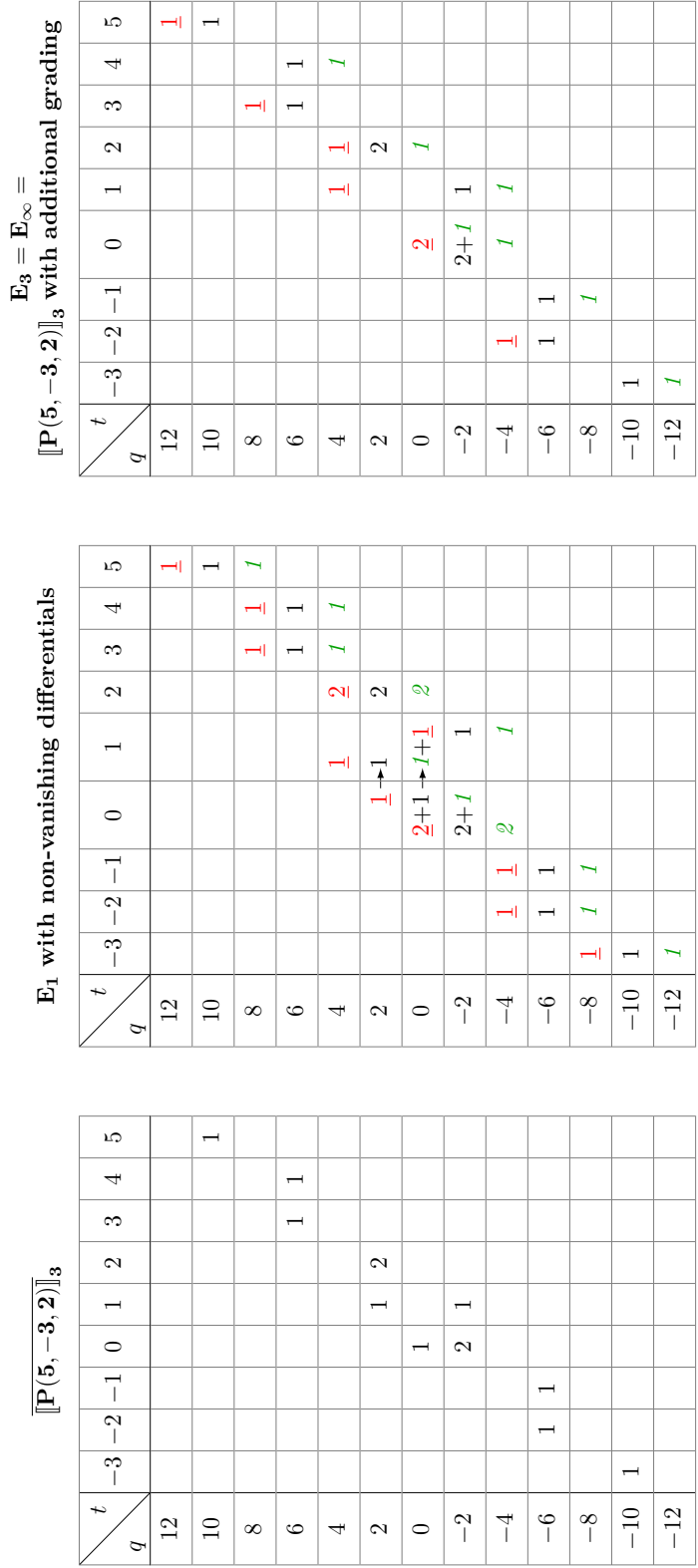
\begin{figure}[p]%
\newcommand{\RH}[1]{\makebox[0.7em]{#1}}%
\rotatebox{90}{\scalebox{0.91}{\parbox{1.1\textheight}{\centering
\extrarowsep=.57ex
\begin{tabu}{|[gray]c|*{9}{c|[gray]}}%
\multicolumn{10}{c}{\raisebox{0pt}[2ex]{$\mathbf{\Hom[r]{P(5,-3,2)}{3}}$}}\medskip \\
\tabucline[gray]-%
\backslashbox{$q$}{$t$} & \RH{$-3$}  & \RH{$-2$}  & \RH{$-1$}  & \RH{0} & \RH{1} & \RH{2} & \RH{3} & \RH{4} & \RH{5} \\\hline
12  &  &   &   &   &   &   &   &   &   \\\tabucline[gray]-
10  &  &   &   &   &   &   &   &   & 1 \\\tabucline[gray]-
8  &  &   &   &   &   &   &   &   &    \\\tabucline[gray]-
6  &  &   &   &   &   &   & 1 & 1 &    \\\tabucline[gray]-
4  &  &   &   &   &   &   &   &   &    \\\tabucline[gray]-
2 &  &   &   &   & 1 & 2 &   &   &     \\\tabucline[gray]-
0 &  &   &   & 1 &   &   &   &   &     \\\tabucline[gray]-
$-2 $ &  &   &   & 2 & 1 &   &   &   &    \\\tabucline[gray]-
$-4 $ &  &   &   &   &   &   &   &   &    \\\tabucline[gray]-
$-6 $ &  & 1 & 1 &   &   &   &   &   &    \\\tabucline[gray]-
$-8 $ &  &   &   &   &   &   &   &   &    \\\tabucline[gray]-
$-10$  & 1 &   &   &   &   &   &   &   &  \\\tabucline[gray]-
$-12$  &  &   &   &   &   &   &   &   &   \\\tabucline[gray]-
\end{tabu}\hfill
\newcommand{\DD}[1]{\emph{\textcolor{darkgreen}{#1}}}%
\newcommand{\UU}[1]{\textcolor{darkred}{\underline{#1}}}%
\newcommand{\Pfeil}[1]{\makebox[0pt][l]{\:%
\setlength{\unitlength}{5mm}%
\begin{picture}(0,0)\put(0,0.18){\vector(1,0){#1}}\end{picture}}}%
\begin{tabu}{|[gray]c|*{9}{c|[gray]}}%
\multicolumn{10}{c}{\raisebox{0pt}[2ex]{$\mathbf{E_1}$ \textbf{with non-vanishing differentials}}}\medskip \\
\tabucline[gray]-%
\backslashbox{$q$}{$t$} & \RH{$-3$}  & \RH{$-2$}  & \RH{$-1$}  & \RH{0} & \RH{1} & \RH{2} & \RH{3} & \RH{4} & \RH{5} \\\hline
12          &        &        &        &            &               &        &        &        & \UU{1} \\\tabucline[gray]-
10          &        &        &        &            &               &        &        &        & 1      \\\tabucline[gray]-
8           &        &        &        &            &               &        & \UU{1} & \UU{1} & \DD{1} \\\tabucline[gray]-
6           &        &        &        &            &               &        & 1      & 1      &        \\\tabucline[gray]-
4           &        &        &        &            & \UU{1}\phantom{+1} & \UU{2} & \DD{1} & \DD{1} &        \\\tabucline[gray]-
2           &        &        &        & \phantom{2+}\UU{1}\Pfeil{0.6}    & 1\phantom{+1} & 2      &        &        &  \\\tabucline[gray]-
0           &        &        &        & \UU{2}+1\Pfeil{0.6} & \DD{1}+\UU{1} & \DD{2} &        &        &        \\\tabucline[gray]-
$-2 $       &        &        &        & 2+\DD{1} & \phantom{+1}1             &        &        &        &        \\\tabucline[gray]-
$-4 $       &        & \UU{1} & \UU{1} & \DD{2}\phantom{+1} & \phantom{+1}\DD{1}        &        &        &        &        \\\tabucline[gray]-
$-6 $       &        & 1      & 1      &            &               &        &        &        &        \\\tabucline[gray]-
$-8 $       & \UU{1} & \DD{1} & \DD{1} &            &               &        &        &        &        \\\tabucline[gray]-
$-10$       & 1      &        &        &            &               &        &        &        &        \\\tabucline[gray]-
$-12$       & \DD{1} &        &        &            &               &        &        &        &        \\\tabucline[gray]-
\end{tabu}\hfill
\begin{tabu}{|[gray]c|*{9}{c|[gray]}}%
\multicolumn{10}{c}{\raisebox{1ex}[2ex][0pt]{\makebox[0pt]{\parbox{6.4cm}{\hfill
$\mathbf{E_3 = E_{\infty}=}$\hfill\mbox{}\\$\mathbf{\Hom{P(5,-3,2)}{3}}$ \textbf{with additional grading}}}}}\medskip \\
\tabucline[gray]-
\backslashbox{$q$}{$t$} & \RH{$-3$}  & \RH{$-2$}  & \RH{$-1$}  & \RH{0} & \RH{1} & \RH{2} & \RH{3} & \RH{4} & \RH{5} \\\hline
12    &&  &  &  &  &  &  &  & \UU{1}  \\\tabucline[gray]-
10    &  &  &  &  &  &  &  &  & 1  \\\tabucline[gray]-
8     &&  &  &  &  &  & \UU{1} &  &      \\\tabucline[gray]-
6     &&  &  &  &  &  & 1 & 1 &     \\\tabucline[gray]-
4     &&  &  &  & \UU{1} & \UU{1} &  & \DD{1} &     \\\tabucline[gray]-
2     &&  &  &  &  & 2 &  &  &       \\\tabucline[gray]-
0     &&  &  & \UU{2} &  & \DD{1} &  &  &      \\\tabucline[gray]-
$-2 $ &&  &  & 2+\DD{1} & 1 &  &  &  &      \\\tabucline[gray]-
$-4 $ && \UU{1} &  & \DD{1} & \DD{1} &  &  &  &     \\\tabucline[gray]-
$-6 $ && 1 & 1 &  &  &  &  &  &      \\\tabucline[gray]-
$-8 $ &&   & \DD{1} &  &  &  &  &  &  \\\tabucline[gray]-
$-10$ & 1 &  &  &  &  &  &  &  &  \\\tabucline[gray]-
$-12$ & \DD{1} &  &  &  &  &  &  &  &  \\\tabucline[gray]-
\end{tabu}\bigskip\\
\protect\caption{%
On the left, the reduced $\mathfrak{sl}_3$-homology of the $P(5,-3,2)$-pretzel knot.
In the middle, the first page of the reduced-unreduced spectral sequence. The $r$-grading
is expressed by colours (and underlining/italics). On this page, ``red may kill black and black green''. On the
second page (not shown), only ``red may kill green''. The total differential on the first and
second page has rank 2 and 5, respectively.
On the right, the third page of homology, which is the limit of the spectral sequence.}%
\protect\label{fig:page}}%
}}%
\end{figure}
This section is devoted to the proof of \thmref{thm:sseq}.
We need the following technical lemma, whose proof is left to the reader.
\begin{Lem}
\label{prop:gradfilt}
Let $(C,\partial)$ be a filtered chain complex, whose filtration we denote by $\mathcal{F}$.
Let there be an additional grading $C = \bigoplus_{i\in\mathbb{Z}} C_i$
that is respected by the differential.
The filtration $\mathcal{F}$ induces a filtration on each $C_i$ by $\mathcal{F}_jC_i := C_i \cap \mathcal{F}_jC$.
If $C$, as a filtered vector space, is the sum of the filtered $C_i$, we say that the
filtration is compatible with the additional grading.
In this case, the spectral sequence induced by $\mathcal{F}$ respects the additional grading on $C$.
\end{Lem}
\begin{Pf}[Proof of {\thmref{thm:sseq}}]
\textbf{(i)}:
Recall that $C_N(D)$ is a module over $A = \mathbb{C}[X]/(X^N)$.
Let us introduce a filtration $\mathcal{R}$ on $A$, given by
$\mathcal{R}_{2i-N+1}A = (X^i)$. Explicitly, we have \[
A = \mathcal{R}_{-N+1} A \supset \mathcal{R}_{-N+3} \supset \ldots \supset \mathcal{R}_{N-1} A \supset \{0\}.
\]
This induces a filtration on $C_N(D)$. Let us denote it by $\mathcal{R}$ as well,
and call the induced grading the $r$-grading.
Since the differential of $C_N(D)$ commutes with the $A$-scalar multiplication,
it also respects this filtration.
So $\mathcal{R}$ induces a spectral sequence $E_{\bullet}$, which
(forgetting the additional degree) converges to $\Hom{L}{N}$.
Note moreover that $\mathcal{R}$ and the $q$-grading are compatible
in the sense of \thmref{prop:gradfilt}.
Hence $E_{\bullet}$ respects the $q$-degree, and its differential on the $k$-th page has
degree $tr^{2k}q^0$.

\textbf{(ii)}:
Let us analyse the $0$-th page of that spectral sequence, i.e. the associated graded chain complex.
We have $\overline{C_N(D)} = C_N(D) \otimes_A \widetilde{\mathbb{C}}$. This is isomorphic to
$C_N(D) / ((X)\cdot C_N(D)) = \mathcal{R}_{-N+1} C_N(D) / \mathcal{R}_{-N+3} C_N(D)$,
with a $q$-shift of $N-1$.
Note that for $i \in \{0, \ldots, N-1\}$,
the multiplication by $X^i$ is an isomorphism 
between $\mathcal{R}_{-N+1} A / \mathcal{R}_{-N+3} A$
and $\mathcal{R}_{2i-N+1} A / \mathcal{R}_{2i-N+3} A$.
Since $C_N(D)$ is a \emph{free} $A$-module, this is true for $C_N(D)$ as well:
the multiplication by $X^i$ is an
isomorphism
\[
\mathcal{R}_{-N+1} C_N(D) / \mathcal{R}_{-N+3} C_N(D) \to
\mathcal{R}_{2i-N+1} C_N(D) / \mathcal{R}_{2i-N+3} C_N(D)
\]
 of complex vector spaces.
Because the $A$-scalar multiplication
commutes with the differential, this map is an isomorphism of chain complexes. It
shifts the $q$-grading and the $r$-grading by $2i$. The $0$-th page of the spectral sequence is the sum
of the $\mathcal{R}_{2i-N+1} C_N(D) / \mathcal{R}_{2i-N+3} C_N(D)$, and thus
\begin{equation*}
E_0 \cong
\bigoplus_%
{i\in\{-N+1,-N+3,\ldots,N-1\}}
(qr)^i\cdot\Red{C}_N(D).
\end{equation*}
Taking homology yields the stated result for the first page.

\textbf{(iii)}:
To prove the invariance of the higher pages, let us use the following
lemma proved e.g. in \cite[theorem 3.5]{usersguide} (also used by Rasmussen \cite[lemma 6.1]{rasmussen}):
\begin{Lem}\label{lem:lowiso}
Let $f: C\to C'$ be a map of filtered chain complexes. Let $E_{\bullet}$ and $E_{\bullet}'$
be the respective spectral sequences associated to $C$ and $C'$, 
and for all $r \geq 0$, let $f_r$ be the induced graded map from $E_r$ to $E_r'$.
If $f_R$ is an isomorphism for some $R$, then $f_r$ is also an isomorphism for all $\infty \geq r \geq R$.
\end{Lem}
So the map of \thmref{prop:red} induces an isomorphism between the higher pages of
the spectral sequence associated to diagrams related by a Reidemeister move.
\end{Pf}
\Thmref{fig:page} shows the reduced-unreduced spectral sequence for the pretzel knot $P(5,-3,2)$.

%% file: tool.tex
\begin{figure}
\newcommand{\littlebox}[1]{%
\fbox{\parbox{3.3cm}{\centering #1\\}}}%
\hspace{-4.75em}%
\xymatrix@C=3.3cm@R=1.2cm{
\littlebox{Reduced \homf\-homology~\Hom[r]{L}{\infty}}
\ar@{-->}[ddd]_-{\parbox{3.2cm}{\raggedleft \Thmref{cor:tool}\ \ }}
\ar@{=>}[r]_-{\parbox{3.2cm}{\centering \footnotesize \raisebox{-2.5ex}{Trivial for large $N$}}}^-{%
\parbox{4cm}{\footnotesize \centering \raisebox{1ex}{$\xdeg d_k = tq^{2Nk}a^{-2k}$}}}
&
\ar[r]_-{\mbox{$a\mapsto q^N$}}
&
\littlebox{Reduced \Sln-homology \Hom[r]{L}{N}}
\ar@{->}[d]^-{\raisebox{0.7ex}{\ \mbox{$\cdot\,[N]_{qr}$}}}
\\
&&
\ar@{=>}[d]_-{\parbox{2.2cm}{\raggedleft $\xdeg d_k = tr^{2k}$}}
\\
&&
\ar@{->}[d]^-{\raisebox{0.7ex}{\ \mbox{$r\mapsto 1$}}}
\\
\littlebox{Filtered \Sln-homology $\Hom{L}{N}^f$}
\ar@{<=}[rr]^-{\raisebox{1ex}{\parbox{4cm}{\centering $\xdeg d_k = t q^{2Nk}$}}}_-{\raisebox{-2ex}{\parbox{4cm}{\centering Conjecture: $E_2 = E_{\infty}$}}}
&&
\littlebox{Unreduced \Sln-homology \Hom{L}{N}}
}
\caption{The proof of \thmref{thm:main1b} and \thmref{cor:tool} in a nutshell. Double arrows stand for spectral sequences, single arrows for other relations.}
\label{fig:schaubild}
\end{figure}
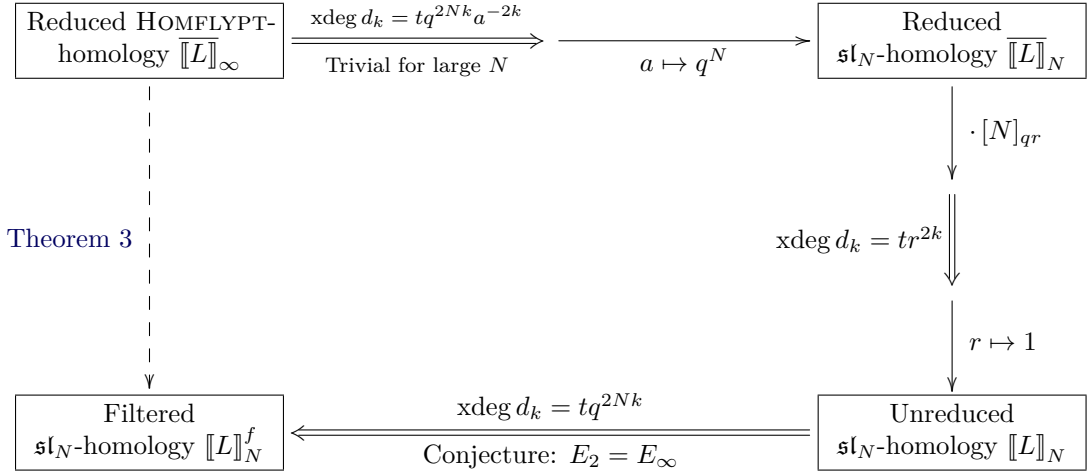

The following lemma (a direct generalisation of \cite[theorem 5.1]{manandmachine})
describes the decategorification of a spectral sequence:
\begin{Lem}\label{lem:sseq}
Let $(E_{\bullet}, d_{\bullet})$ be a spectral sequence of $\Int^n$-graded finite dimensional vector spaces.
Then for all $k\geq 1$ there are polynomials $f_k \in \Nat[x_1^{\pm 1},\ldots x_n^{\pm 1}]$, such that
for all $\ell\geq 1$ the following decomposition holds:
\begin{equation*}
\xdim E_{\ell} = \xdim E_{\ell+1} + \sum_{k=1}^{\ell} (1 + \xdeg d_k)\cdot f_k.
\end{equation*}
In particular,
\begin{equation*}
\xdim E_1 = \xdim E_{\infty} + \sum_{k=1}^{\infty} (1 + \xdeg d_k)\cdot f_k.
\end{equation*}
The spectral sequence converges on the $\ell$-th page if and only if
$\forall k \geq \ell: f_k = 0$.
\end{Lem}
\begin{Prop}\label{thm:main1b}
Let $K$ be a knot, and let $N \geq 2$.
There are polynomials $\Red{P}'_N \in \Nat[t^{\pm 1},q^{\pm 1},a^{\pm 1}]$, $P'_N\in \Nat[t^{\pm 1},q^{\pm 1},r^{\pm 1}]$ and for all
$k \geq 1$ polynomials $f_N^k \in \Nat[t^{\pm 1},q^{\pm 1},a^{\pm 1}]$, $g_N^k\in\Nat[t^{\pm 1},q^{\pm 1},r^{\pm 1}]$, $h_N^k\in\Nat[q^{\pm 1},t^{\pm 1}]$,
such that for large enough $N$ we have $\forall k: f_N^k = 0$, and
such that the following decompositions hold:
\begin{align*}
\xdim \Hom[r]{K}{\infty} & = \Red{P}'_N(t,q,a) + \sum_{k=1}^\infty (1 + tq^{-2Nk}a^{2k})f_N^k(t,q,a), \\
\xdim \Hom[r]{K}{N} & = \Red{P}'_N(t,q,q^N), \\
\xdim \Hom[r]{K}{N} \cdot [N]_{qr} & = P'_N(t,q,r) + \sum_{k=1}^{N-1} (1 + tr^{2k})g_N^k(t,q,r), \\
\xdim \Hom{K}{N} & = P'_N(t,q,1), \\
\xdim \Hom{K}{N} & = q^{s'_N(K)}\cdot [N]_q + \sum_{k=1}^{\infty} (1 + tq^{2Nk})h_N^k(q,t).
\end{align*}
\end{Prop}
\begin{Pf}
Apply \thmref{lem:sseq} to Rasmussen's spectral sequence (\thmref{prop:Rasmussens-spectral-sequence}),
the reduced-unreduced spectral sequence (\thmref{thm:sseq}) and the Lee-Gornik spectral sequence with
renumbered pages (\thmref{prop:Lees-spectral-sequence}). See also \thmref{fig:schaubild}.
\end{Pf}
Given $\xdim\Hom[r]{K}{\infty}$, there are only finitely many choices for the auxiliary polynomials in the proposition
and for the $s'_N(K)$. So \homf-homology induces restrictions on the \Sln-concordance invariants (and thus a lower bound on the slice genus).
\Thmref{thm:main1b} may appear unwieldy; and in fact, we will only use \thmref{cor:tool}, which skips all
intermediary steps between the \homf-homology and the \Sln-concordance invariants.
\begin{Pf}[Proof of {\thmref{cor:tool}}]
Let us use the equations of \thmref{thm:main1b}, climbing from the bottom up:
\begin{align*}
q^{s_N'(K)-N+1} & \leq \xdim \Hom{K}{N}^f  \\
\displaybreak[0]
\implies\quad q^{s_N'(K)-N+1} & \leq \xdim \Hom{K}{N}  \\
\displaybreak[0]
\implies\quad q^{s_N'(K)-N+1}r^i & \leq P_N' \text{ for some $i$}  \\
\displaybreak[0]
\implies\quad q^{s_N'(K)-N+1}r^i & \leq \xdim\Hom[r]{K}{N} \cdot [N]_{qr}  \\
\displaybreak[0]
\implies\quad q^{s_N'(K)-N+1+j} & \leq \xdim\Hom[r]{K}{N} \text{ for some $j\in\{1-N, 3-N, \ldots, N-1\}$}  \\
\displaybreak[0]
\implies\quad q^{\alpha}a^{\beta} & \leq \Red{P}'_N \text{ for some $\alpha, \beta$ with $\alpha + N\beta \leq s_N'(K)$ }  \\
\displaybreak[0]
\implies\quad q^{\alpha}a^{\beta} & \leq \xdim\Hom[r]{K}{\infty}.
\end{align*}
An analogous reasoning yields $q^{s'_N(K)+N-1} \leq \xdim\Hom{K}{N} \implies q^{\alpha'}a^{\beta'}\leq \Hom[r]{K}{\infty}$ for some $\alpha', \beta'$ with
$s'_N(K) \leq \alpha' + N\beta'$.
\end{Pf}
Note that the power of \thmref{thm:main1b} is limited:
\begin{Prop}\label{prop:weakness}
Let $K$ be a knot, and let $N \geq 2$. Suppose there are polynomials
$i^0_N, \ldots i^{N-1}_N \in \Nat[t^{\pm 1}, a^{\pm 1}, q^{\pm 1}]$
and some $\alpha, \beta$ such
that the following decomposition holds:
\begin{equation*}
\xdim \Hom[r]{K}{\infty} = q^{\alpha}a^{\beta} + \sum_{k = 0}^{N-1} (1 + ta^{-2}q^{2k})\cdot i^k_N.
\end{equation*}
Then there is also a decomposition as in \thmref{thm:main1b} with $\alpha + N\beta$
at the place of $s'_N(K)$, i.e. the theorem cannot be used to show $s'_N(K) \neq \alpha + N\beta$.
\end{Prop}
\begin{Pf}
Let $N$ be fixed. Setting $f^k_N = 0$ for all $k$, $h^k_N = 0$ for $k \geq 2$ and
\begin{align*}
h^1_N & = \sum_{k=0}^{N-1} [N-k]_q \cdot q^k\cdot i^k_N, \\[\medskipamount]
g^k_N & = [N-k]_{qr} \cdot (qr)^{-k}\cdot i^{N-k}_N \\
\end{align*}
gives the desired decomposition.
\end{Pf}
\begin{Rmk}
If Rasmussen's spectral sequence from $\Hom[r]{K}{\infty}$ to the regraded version of $\Hom[r]{K}{1}$
converges on the second page, it gives a decomposition as in the above proposition,
with $i^k_N = 0$ for $k \neq 1$ and $\alpha = - \beta$.
So for any knot for which that spectral sequence converges on the second page,
\thmref{thm:main1b} alone is not strong enough to distinguish
the \Sln-concordance invariants.
\end{Rmk}

%% file: slicetorus.tex
This section is largely independent from sections 2--4. It collects and extends what is known about slice-torus invariants.
Let us start by listing some invariants that are not slice-torus:
the classical knot signature $\sigma$ or the concordance invariant $\delta$
from the Floer homology of double branched covers \cite{delta}
give slice genus bounds, but they are not sharp for torus knots;
the Lipshitz-Sarkar invariant \cite{lipsar} on the other hand is not a concordance homomorphism.

The following results are mainly due to Livingston \cite{livingston}, who built on the work of Rudolph \cite{rudolph}.
Although Livingston only considers slice-torus invariants which take even integer
values, his results and proofs carry over unchanged to the general case of real-valued invariants.
Note that we opted for a different normalisation of the slice-torus invariants.
Throughout this section, let $\nu$ denote an arbitrary slice-torus invariant.
The proof of the following proposition is standard:
\begin{Prop}\label{prop:cob}
\begin{renumerate}
\item[\textbf{(i)}] ({{\cite[Cor. 2]{livingston}}})
For all knots $K$, the absolute value of $\nu$ is a lower bound to twice the slice genus: $|\nu(K)| \leq 2g_4(K)$.
\item[\textbf{(ii)}]
If there is a connected smooth cobordism of Euler characteristic $\chi$ between two knots $K_0$ and $K_1$,
then
\begin{equation*}
|\nu(K_0) - \nu(K_1)| \leq -\chi.
\end{equation*}
\end{renumerate}
\end{Prop}
\begin{figure}[h]
\centering
\avcfig{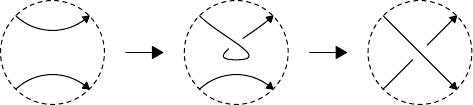}
\caption[A cobordism of Euler characteristic $-1$ inserting a positive crossing.]%
{A cobordism of Euler characteristic $-1$ inserting a positive crossing.
The cobordism consists of a Reidemeister I move and a saddle move.}
\label{fig:crossing cob}
\end{figure}
\begin{Lem}\label{lem:basiccobo}
There is %
a cobordism of Euler characteristic $-1$ inserting or resolving a positive or negative crossing,
\end{Lem}
\begin{Pf}
See \thmref{fig:crossing cob}.
\end{Pf}
\begin{Prop}[{{\cite[Cor. 3]{livingston}}}]\label{thm:crossing-change}
If $K_+$ and $K_-$ are knots that have diagrams that are identical but
for the sign of one crossing, which is given in the subscript (see \thmref{fig:Skein-Homf2}), then
\[
0 \leq \nu(K_+) - \nu(K_-) \leq 2.
\]
\end{Prop}
\begin{Lem}[{\cite[Cor. 7]{livingston}}]
\label{cor:posBraid}
Let $B$ be a positive braid, i.e. a braid whose word contains only the $\sigma_i$,
not the $\sigma_i^{-1}$. Suppose $B$ has $n$ strands and $k$ crossings.
If the closure $\tr(B)$ of $B$ is a knot, then $\nu(\tr(B)) = 2g_4(\tr(B)) = 2g_3(\tr(B)) = 1 + k - n$.
\end{Lem}
Quasi-positivity has been introduced and studied by Rudolph, see \cite{rudolph}.
\begin{Def}
A braid $B$ is said to be \emph{quasi-positive} if it is the product of braid-words that
are conjugate to one of the $\sigma_i$; i.e. $B = \prod_j w_j \sigma_{i_j} w_j^{-1}$,
where $w$ is any braid-word.
\end{Def}
\begin{Prop}[\cite{rudolph}]
Let $B$ be a quasi-positive braid with $k_+$ positive crossings,
$k_-$ negative crossings, writhe $w = k_+ - k_-$ and $n$ strands.
If $\tr(B)$ is a knot,
then $2g_4(\tr(B)) = 1 + w - n$.
\end{Prop}
The following \namecref{thm:qpos} has been proven for the Rasmussen invariant by Shumakovitch \cite{shumakovitch};
for $2\tau$ it is an immediate consequence of the results of Plamenevskaya \cite{plamenevskaya}.
The relationship between the $\tau$-invariant, quasi-positivity and fibredness were studied by Hedden \cite{hedden_positivity}.
\begin{Prop}\label{thm:qpos}
Slice-torus invariants detect the slice-genus of quasi-positive knots.
\end{Prop}
\begin{Pf}
Let $B$ be a quasi-positive braid.
Let $B'$ be the braid obtained by switching every negative crossing of $B$ to a positive one.
By \thmref{cor:posBraid}, $\nu(\tr(B')) = 1 + k_+ + k_- -n$, and
by \thmref{thm:crossing-change},
$\nu(\tr(B')) - \nu(\tr(B)) \leq 2k_- \implies \nu(\tr(B)) \geq 1 + w - n$.
But $1 + w - n = g_4(\tr(B))$ by the previous proposition.
\end{Pf}
\begin{Cor}\label{lem:poslink}
Let $D$ be a positive knot diagram of a knot $K$, i.e. a diagram with only positive crossings.
Then $\nu$ detects the slice genus of $K$.
\end{Cor}
\begin{Pf}
Follows from the previous proposition since positive links are quasi-positive \cite{rudolph2}.
\end{Pf}
\begin{Prop}[{\cite{rasmussen}}]
Let $D$ be a positive knot diagram of a knot $K$, with $k$ crossings, and $n$ Seifert circles.
Then $g_4(K) = 1 + k - n$.
\end{Prop}

One of the strongest restrictions that can be deduced from the slice-torus conditions
is an inequality à la Bennequin \cite{bennequin}.
The first version was proven by Rudolph \cite{rudolph} for the slice-genus,
and by Rasmussen \cite{rasmussen},
Shumakovitch \cite{shumakovitch} and Plamenevskaya \cite{plamenevskaya} for the Rasmussen invariant.
It was subsequently sharpened by Kawamura \cite{kawamura}; her version was generalised by Wu to the $s_N$-invariants \cite{wu3}.
Then, it was  honed yet more independently by Lobb \cite{lobbineq} and Kawamura \cite{kawamuraNotPub}.
Given a diagram $D$ of a knot $K$, the sharper slice-Bennequin inequality gives an upper and lower bound for $\nu(K)$.
\newcommand{\SG}{\Gamma}
\newcommand{\CL}{\tr}
Those bounds are easily computable from $D$, depending only the \emph{Seifert graph} $\SG(D)$:
\begin{Def}%
\label{def:sg}%
The \emph{Seifert graph} $\SG(D)$ of a link diagram $D$ is a planar bipartite graph whose edges carry a sign ($+$ or $-$).
It is constructed as follows:
the vertices of $\SG(D)$ correspond to the circles of the Seifert resolution of $D$.
A fixed crossing of $D$ is adjacent to two different Seifert circles, which correspond to two
vertices in $\SG(D)$. For any crossing, let $\SG(D)$ have an edge between these two vertices.
The edge's sign indicates if the crossing is positive or negative.
Let $\SG^+(D)$ (\,$\SG^-(D)$) be the subgraph of $\SG(D)$ that contains only the positive (negative) edges.
Let $O^{\pm}(D)$ be the number of connected components of $\SG^{\pm}(D)$.
\end{Def}
\begin{Thm}[The sharper slice-Bennequin inequality]\label{lem:lobbineq}
Let $D$ be a diagram of a knot $K$, with writhe $w$ and $n$ Seifert circles.
Then
\begin{equation*}
- 1 + w - n + 2O^+
\leq \nu(K) \leq
1 + w + n - 2O^-.
\end{equation*}
\end{Thm}
\begin{figure}
\newcolumntype{x}[1]{%
>{\centering\hspace{0pt}}p{#1}}%
\begin{tabular}{*{4}{x{0.22\textwidth}}}
\lowerhalfx{\reflectbox{\includegraphics[height=3.7cm]{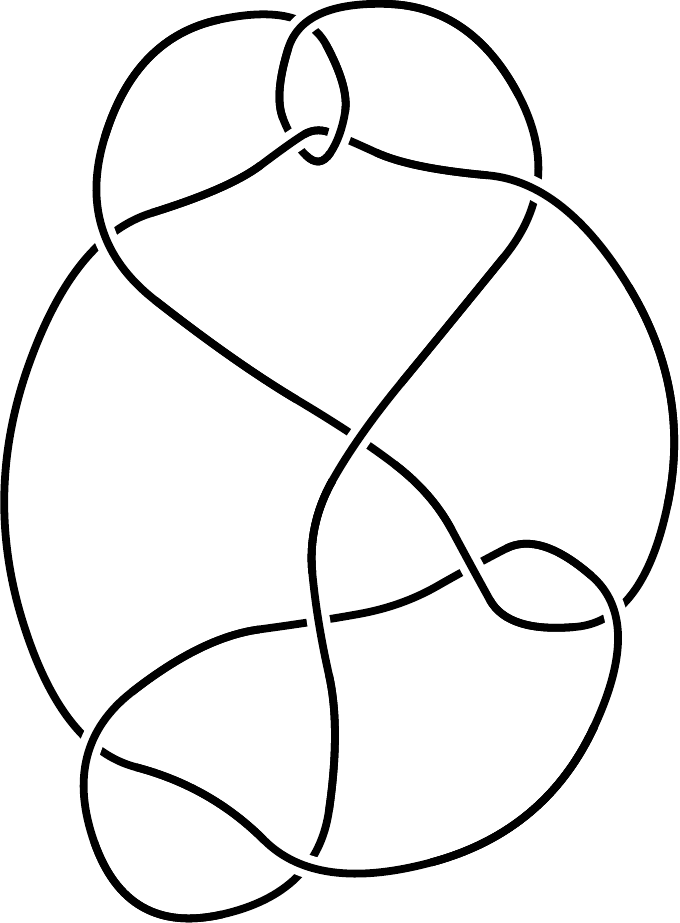}}} &
\lowerhalfx{\reflectbox{\includegraphics[height=3.7cm]{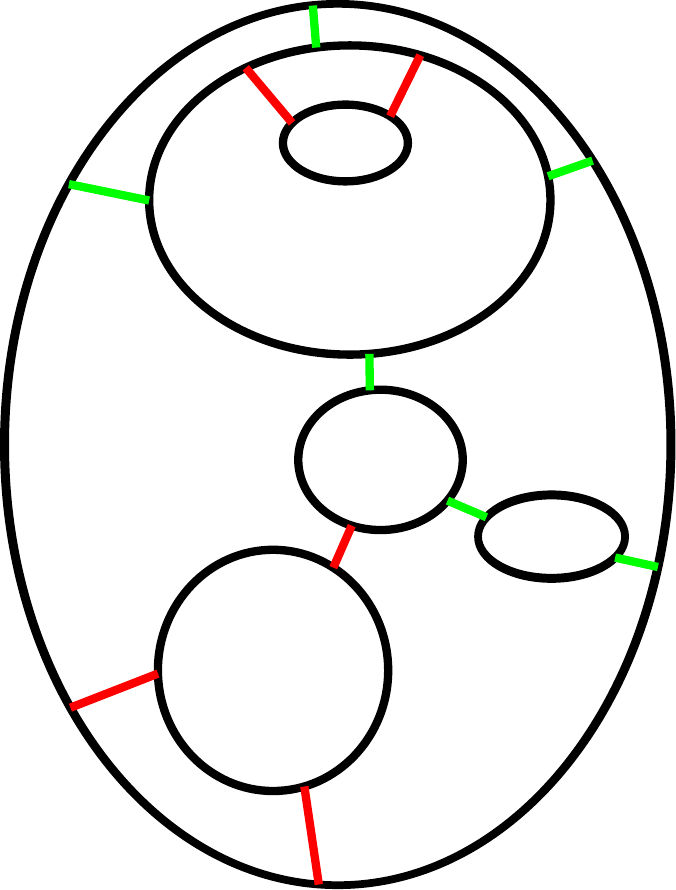}}} &
\lowerhalfx{{\includegraphics[height=3.7cm]{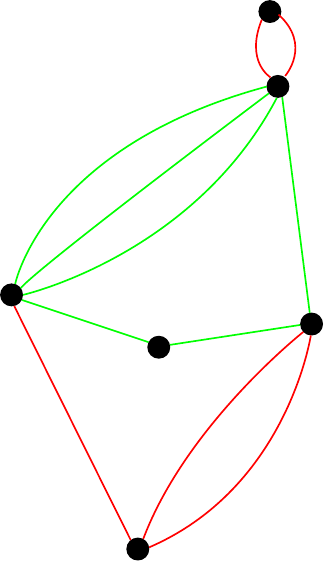}}} &
\lowerhalfx{{\includegraphics[height=3.7cm]{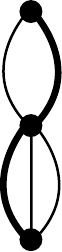}}} \tabularnewline
\rule{0ex}{3.5ex} Diagram $D$ of $K$ & Seifert resolution & $\Gamma(D)$ &  $G(D)$, $T$ drawn bold
\end{tabular}\\[\baselineskip]
\begin{tabular}{*{3}{x{0.28\textwidth}}}
\lowerhalfx{\reflectbox{\includegraphics[height=3.7cm]{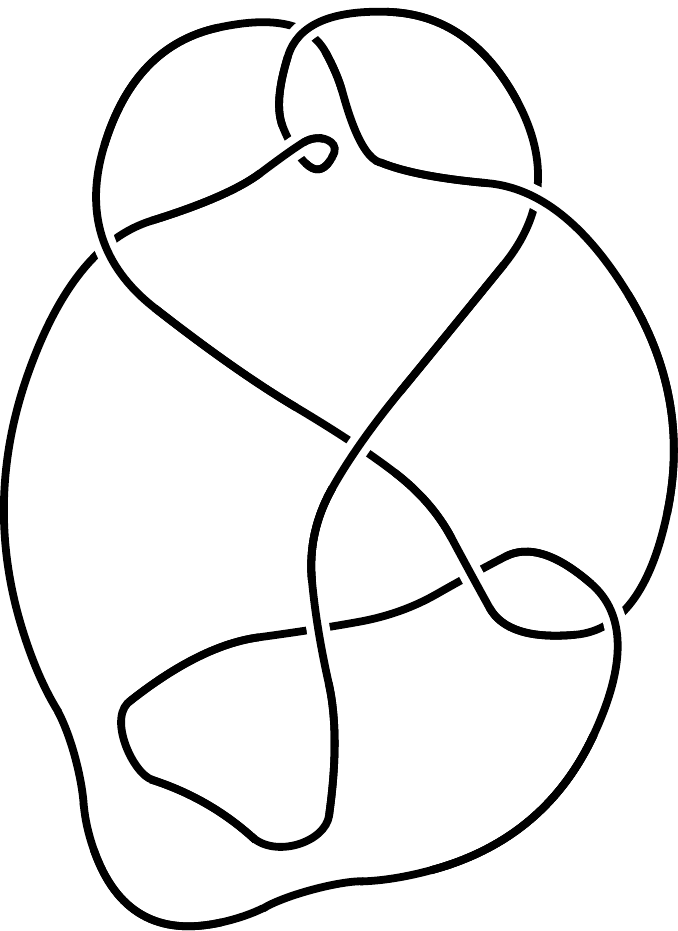}}} &
\lowerhalfx{\reflectbox{\includegraphics[height=3.7cm]{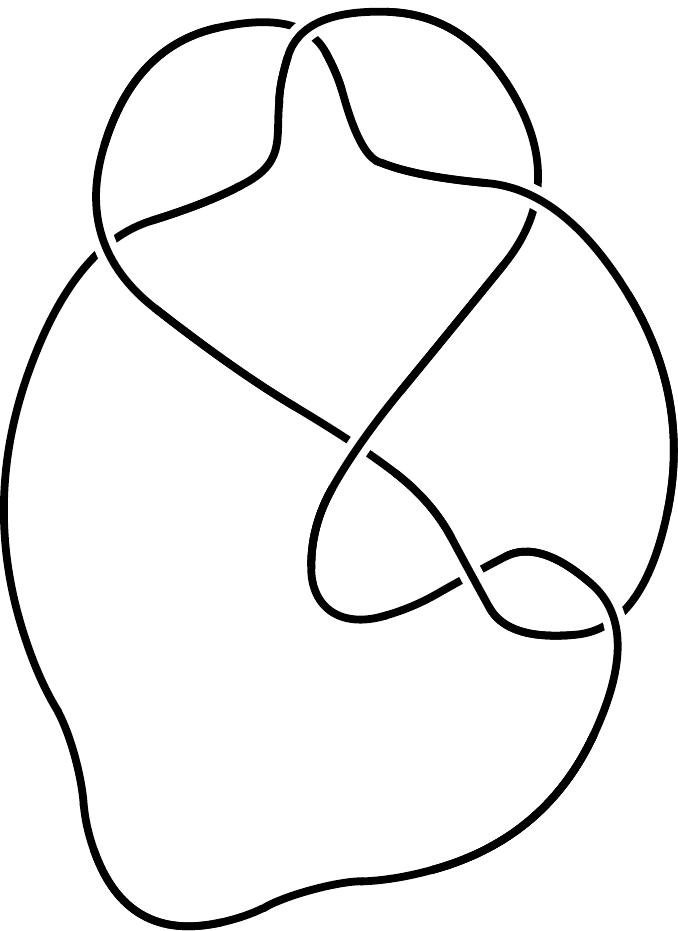}}} &
\lowerhalfx{\reflectbox{\includegraphics[height=3.7cm]{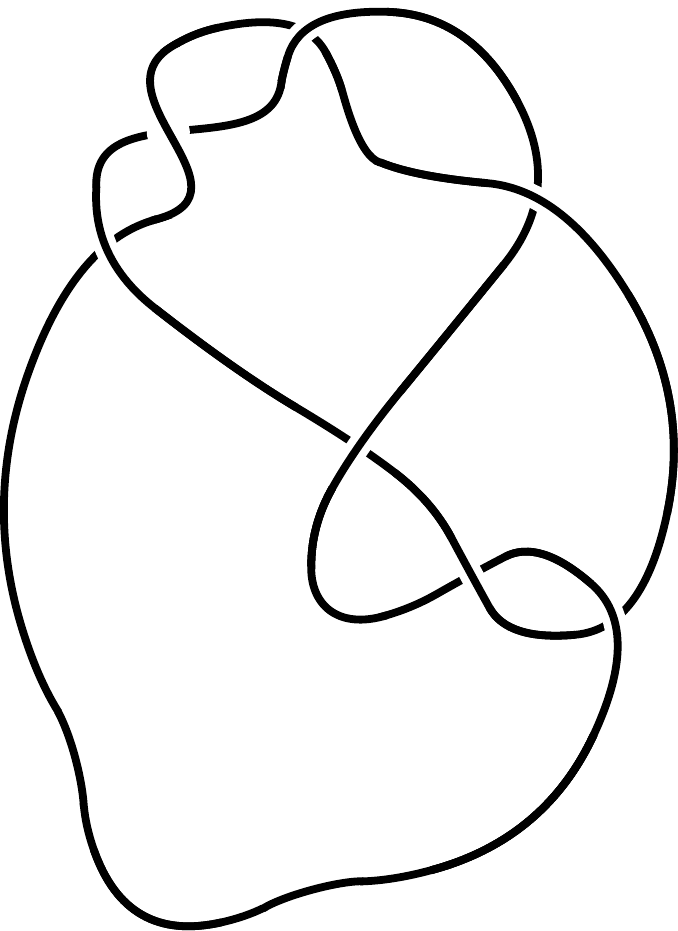}}}\tabularnewline
\rule{0ex}{3.5ex}$D'$ & $D''$ & Diagram $D'''$ of $K'$
\end{tabular}
\caption{%
A diagram of the knot $11n_{53}$, and what happens to it in the proof of \thmref{lem:lobbineq}.
Diagrams drawn with knotscape \cite{knotscape}.}
\label{fig:53}
\end{figure}
\begin{Pf}[Proof (following \cite{abe}).]
Let us only prove the lower bound, since the upper bound then follows from $\nu(-K) = -\nu(K)$.
For an example of the following constructions, see \thmref{fig:53}.
Let $G(D)$ be the graph obtained from $\SG(D)$ by contracting all positive edges; or more explicitly, the graph
that has as vertices the components of $\SG^+(D)$, and has for each negative
edge in $\SG(D)$ an edge between the corresponding vertices of $\SG^+(D)$. Then $G(D)$ is a connected
graph with $O^+$ many vertices and $k^-$ many edges. Pick $O^+ - 1$ edges that form a tree $T$.
Gluing together $k^- - O^+ + 1$ many copies of the cobordism of \thmref{lem:basiccobo}
gives a connected cobordism with Euler characteristic $- k^- + O^+ - 1$ between $K$ and a link with diagram $D'$,
such that $G(D')$ is the tree $T$.
The $O^+ - 1$ many negative crossings of $D'$ are thus nugatory and may be
removed by twists; each twist diminishes the number of Seifert circles by one. The
ensuing diagram $D''$ is positive, with $k^+$ many crossings and $n - O^+ +1$ many Seifert circles.
If the link represented by $D''$ has $c$ components, connect them by adding $c - 1$ positive crossings using the
cobordism of \thmref{lem:basiccobo}. This gives a cobordism to a knot $K'$ with diagram $D'''$ of
Euler characteristic $1 - c$. Overall, there is a cobordism between $K$ and $K'$ of
Euler characteristic $ - c - k^- + O^+$.
The diagram $D'''$ is positive with $k^+ + c - 1$ many crossings and $n - O^+ + 1$ many Seifert circles,
hence $\nu(K') = k^+ + c - n + O^+ - 1$. 
Finally, by \thmref{prop:cob}(ii),
\begin{align*}
\nu(K') - \nu(K)
& \leq c + k^- - O^+ \\
\implies \phantom{\nu(K') - }\nu(K) & \geq (k^+ + c - n + O^+ - 1) - c - k^- + O^+ \\
\implies \phantom{\nu(K') - }\nu(K) & \geq 1 + w - n + 2O^+ \\[-2\baselineskip]
\end{align*}%
\end{Pf}%
\begin{Cor}\label{cor:alt}
On an alternating knot, slice-torus invariants take the same value as the signature.
\end{Cor}
\begin{Pf}
Follows from \cite[proposition 3.3]{lee-signature}.
\end{Pf}
If $D$ is a knot diagram with $O^+ + O^- = n + 1$, then the lower bound of \thmref{lem:lobbineq} equals the upper bound,
and thus $\nu$ is determined by the inequalities. Such diagrams are called \emph{homogeneous}, and, consequently, a link is called homogeneous
if it has a homogeneous diagram. This notion was introduced by Cromwell \cite{cromwell-homogeneous} and
its relationship with the Rasmussen invariant was studied by Abe \cite{abe}.

For the sake of completeness, let us cite results stated in \cite{doubled}
and \cite{cott} about $2\mathbb{Z}$-valued slice-torus invariants of certain satellite knots.
The results and their proofs (which are therefore omitted)
carry through mostly unchanged to real-valued slice-torus invariants.
For a knot $K$, let $D_{\pm}(K, t)$ be the $t$-twisted positive or negative Whitehead double.
Notice that $g_4(D_{\pm}(K,t)) \leq 1$, and thus $\nu(D_{\pm}(K, t)) \in [-2, 2]$.
Let $TB(K)$ be the Thurston-Bennequin number. Then:
\begin{Prop}\cite{doubled}
\begin{renumerate}
\item $|\nu(D_-(K, t))| + |\nu(D_+(K, t))| \leq 2$.
In particular, $\nu(D_+(K, t)) = \pm 2\Rightarrow \nu(D_-(K,t)) = 0$.
\item
Let $N: \Int \to \mathbb{R}$ be given by $t \mapsto \nu(D_+(K,t))$.
Then $N$ is non-increasing; for $t \leq TB(K)$, we have $N(t) = 2$,
and for $t \geq -TB(-K)$, we have $N(t) = 0$.
\end{renumerate}
\end{Prop}
For two coprime integers $m,n$, let $K_{m,n}$ be the $(m,n)$-cable of $K$,
i.e. the satellite with companion $K$ and pattern the $(m,n)$-torus knot.
\begin{Prop}\cite{cott}
Let us fix a knot $K$ and some $m > 0$, and define $h: \mathbb{Z}\to\mathbb{R}$
by $h(n) = \nu(K_{m,n}) - (m-1)\cdot n$.
Then $h$ is non-increasing and bounded, and $\sup h - \inf h \leq 2(m-1)$.
\end{Prop}

%% file: examples.tex
\label{sec:calc}
This section contains the proofs of \thmref{cor:examples} and \thmref{c1}.
Pretzel knots are a practical family of candidates to disprove the conjecture that all the \Sln-concordance invariants
are equal: they show sufficiently complex behaviour, yet their diagrams allow easy calculations, because they invite an inductive approach.
\begin{Rmk}
In the proof of \thmref{cor:examples}, we will only use \thmref{thm:sseq} and not the potentially stronger \thmref{thm:main1b}.
But even \thmref{thm:main1b} would not be strong enough to completely determine
the value of $s_N(P(\ell, -m, 2))$. For example,
using Webster's programme \cite{krm2} or the skein long exact sequence, one finds that
\begin{multline*}
\Hom[r]{P(5,-3,2)}{\infty} = 
t^{-3}	a^2	q^{4} + 
t^{-2}		q^{6} + 
t^{-1}	a^2	 + 
(		2q^{2} + 
1)		 + \\
t(	a^{-2}	q^{4} + 
	a^2	q^{-4}) + 
2t^2		q^{-2} + 
t^3	a^{-2}	 + 
t^4		q^{-6} + 
t^5	a^{-2}	q^{-4}.
\end{multline*}
But this polynomial has several different decompositions as in \thmref{prop:weakness},
among them one with $\alpha = \beta = 0$, and one with $\alpha = 2, \beta = 0$.
\end{Rmk}
Let us start by verifying the stated values for the Rasmussen invariant $s_2$.
\begin{Pf}[Proof of {\thmref{cor:examples}(i) \& (iii)}]
Khovanov homology of three stranded pretzels has been completely computed, see \cite{manion}.
Since pretzel knots have homological width 3, their Khovanov homology determines their Rasmussen invariant.
But the essential case $\ell > m > n \geq 2$ has a quicker proof:
in that case, the $(\ell,-m,n)$-pretzel knot is quasi-alternating 
(see Champanerkar and Kofman \cite{quasialtpretzels} and Greene \cite{greene}),
and hence its $s_2$--invariant (and twice its $\tau$-invariant) equals its signature (see \thmref{rmk:quasialt}):
\[
s_2(P(\ell, -m, n)) = 2\tau(P(\ell, -m, n)) = \sigma(P(\ell, -m, n)) = \ell - m.
\]
This value of the signature can be easily computed using G\"oritz matrices and
the formula of Gordon and Litherland \cite{gordon}.
\end{Pf}
Let us continue by calculating the higher \Sln-concordance invariants.
\begin{Lem}\label{lem:torus}
For odd $\ell \geq 3$, we have
\begin{equation*}
\Hom[r]{T(\ell, 2)}{\infty} = a^{1-\ell}q^{\ell-1}\cdot\left(
1 + (t^2q^{-4} + t^3a^{-2}q^{-2})\cdot\frac{t^{\ell - 1}q^{2 - 2\ell} - 1}{t^2q^{-4} - 1}
\right).
\end{equation*}
\end{Lem}
\begin{Pf}
First, one may inductively calculate the \homf-polynomial of $T(\ell, 2)$, using its
defining skein relation.
Then, since the $(\ell, 2)$-torus knot is two-bridge, \thmref{prop:twothin} gives
the \homf-homology.
\end{Pf}
\begin{Lem}\label{lem:special-pretzels}
For all $N \geq 2$ and odd $\ell \geq 5$,
\begin{align*}
s_N(P(\ell, 2 -\ell, 2)) & \in \left\{0, \frac{2}{N-1}\right\} \\[\medskipamount]
s_N(P(\ell, 2 -\ell, 4)) & \in 
\begin{cases}
\{0, 2\} & N = 2, \\
[\frac{4}{N-1} - 2, 0] & N \geq 3.
\end{cases}
\end{align*}
\end{Lem}
\begin{Pf}
Let $K_- = P(\ell, 2-\ell, 2)$.
Switching one of the two negative crossings of the last pretzel strand,
one obtains the sum of two torus knots: $K_+ = T(\ell, 2) \# T(2 - \ell,2)$.
Resolving that crossing, one obtains the positive Hopf link (to get its standard diagram,
apply $(\ell - 2)$ Reidemeister II moves).
The homology of those torus knots has been computed in \thmref{lem:torus}, and
the homology is well-behaved with respect to the connected sum (see \thmref{lem:sum}).
This gives us $\xdim\Hom[r]{K_+}{\infty}$. The totally reduced homology of the Hopf link
is known, too, see \thmref{sec:overview}.
So using the skein long exact sequence (see \thmref{prop:homf les}), one finds that
\[
\xdim\Hom[r]{P(\ell,2-\ell,2)}{\infty}^0 \leq (\ell-2)q^{-2} + 1.
\]
By virtue of \thmref{cor:tool}, this proves the first statement of the lemma.
Notice also that
\[
\xdim\Hom[r]{P(\ell,2-\ell,2)}{\infty}^2 \leq (\ell-4)q^{-6} + a^{-2}.
\]

Now let $K_- = P(\ell, 2 -\ell, 4)$, and fix one of the negative crossings of the last pretzel strand.
Then $K_+ = P(\ell, 2 -\ell, 2)$, and once again $L_0$
is the positive Hopf link. So
\[
\xdim\Hom[r]{P(\ell,2-\ell,4)}{\infty}^0 \leq (\ell-4)a^2q^{-6} + 2.
\]
Applying \thmref{cor:tool} concludes the proof of the second statement.
\end{Pf}
\begin{Lem}\label{lem:lobb for pretzel}
Let $\nu$ be any slice-torus invariant, $\ell$ and $m$ be odd positive integers, and $n$
an even positive integer. Then
\[
\nu(P(\ell, -m, n)) \in [\ell - m - 2, \ell - m].
\]
\end{Lem}
\begin{Pf}
The standard diagram of the $(\ell,-m,n)$-pretzel knot, as shown exemplarily in \thmref{fig:pretzel},
has writhe $(\ell - m - n)$, $(n + 1)$ many Seifert circles, $O^+ = n$ and $O^- = 1$.
So the statement follows from the sharper slice Bennequin inequality (\thmref{lem:lobbineq}).
\end{Pf}
Let us now assemble the proof:
\begin{Pf}[Proof of {\thmref{cor:examples}(ii) \& (iv)}]
Let $\ell > m \geq 3$ be odd and $n \geq 4$ even, and let $N \geq 3$.
By \thmref{lem:special-pretzels}, we have
\[
s_N(P(m + 2, -m, 4)) \in \left[\frac{4}{N-1} -2,0\right].
\]
It takes $\frac{n - 4}{2}$ many crossing switches from positive to negative,
and $\frac{\ell - m - 2}{2}$ many crossing switches from negative to positive
to go from $P(m+2, -m, 4)$ to $P(\ell, -m, n)$. Thus by \thmref{thm:crossing-change} we have
\[
s_N(P(\ell, -m, n)) \in \left[\frac{4}{N-1}+2-n,\ell -m - 2\right].
\]
But by \thmref{lem:lobb for pretzel}, $s_N(P(\ell, -m, n)) \in [\ell - m - 2, \ell - m]$.
This leaves $s_N(P(\ell, -m, n)) = \ell - m - 2$ as only value in the intersection of the
two intervals.

Let us now consider the special case $n = 2$. By the same method one finds that
\[
s_N(P(\ell,-m,2)) \in [\ell - m - 2, \ell - m - 2 + 2/(N-1)].
\]
The intersection of this interval with $\frac{2}{N-1}\mathbb{Z}$ leaves
$\ell - m - 2$ and $\ell - m - 2 + 2/(N-1)$ as only possible values.
\end{Pf}
\begin{figure}
\includegraphics[scale=0.35]{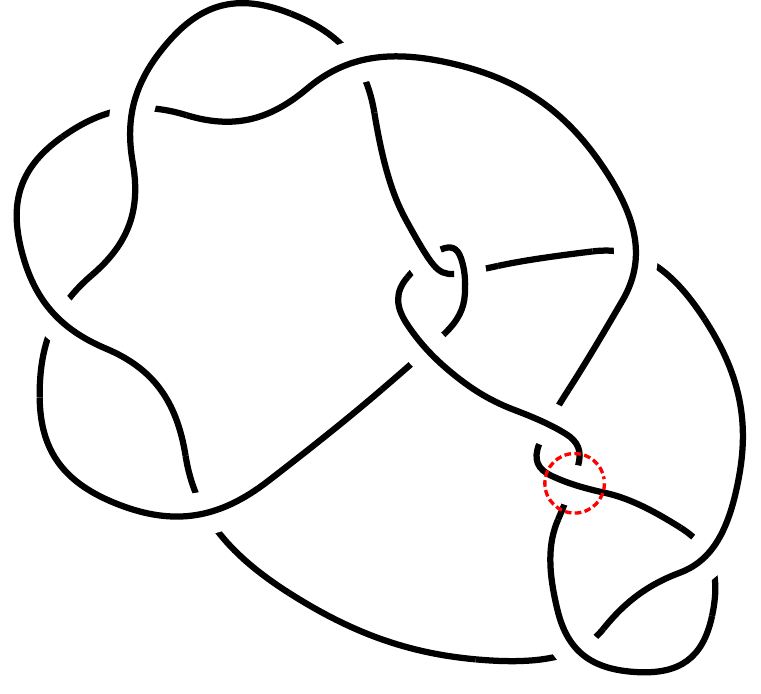}
\hfill
\includegraphics[scale=0.35]{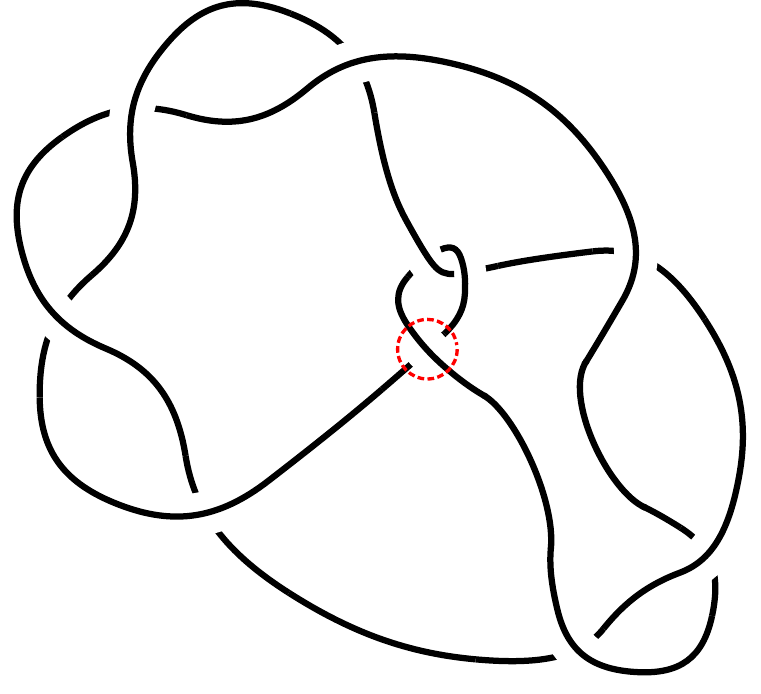}
\hfill
\includegraphics[scale=0.35]{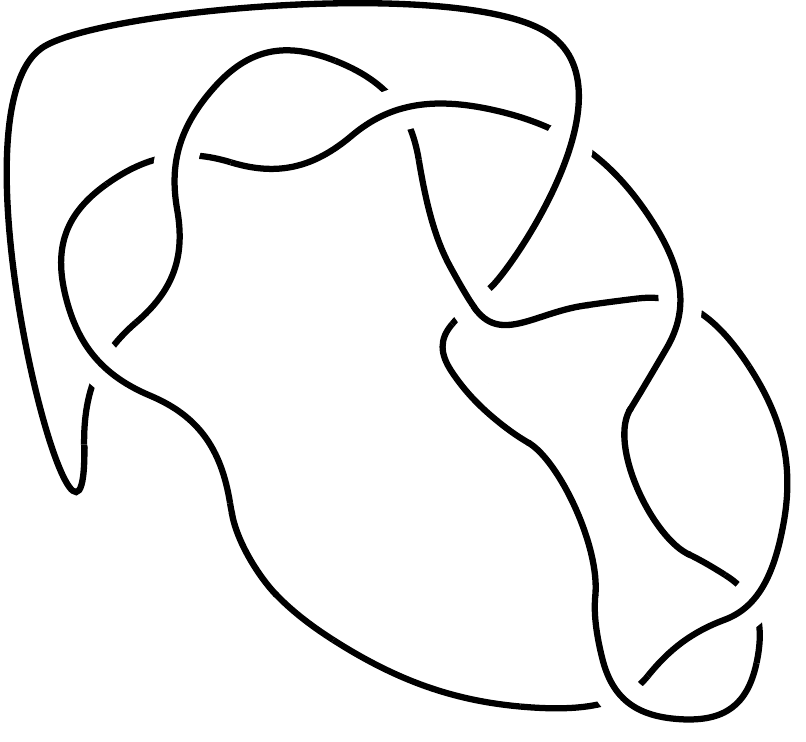}
\caption{From left to right: $12n_{340}$, $10_{141}$ and $8_9$ (drawn with Knotscape \cite{knotscape}).}
\label{fig:12n340}
\end{figure}
The linear independence of $\tau, s_2$ and $s_3$ from the $s_N$ with $N\geq 4$ now follows quickly.
\begin{Pf}[Proof of {\thmref{c1}}]
\textbf{(i)}:
To prove the first statement, note that any linear combination of
$\{s_N\}_{N\geq 3}$ vanishes for $P(7,-5,4)$, but $s_2$ and $\tau$ do not.
Concerning the second statement, invariance of $s_2$ and $\tau$ is due to \cite{sneqtau},
so it suffices to produce a knot $K$ with $s_2(K) = \tau(K) = 0$, but $s_N(K) \neq 0$.
This is accomplished by the quasi-alternating knot $K = P(5,-3,2) \# -T(3,2)$.

\textbf{(ii)}:
Consider a linear combination $u$ of $\{s_N\}_{N\geq 4}$. For $u$ to be equal to $s_3$,
the linear combination needs to be convex (i.e. the sum of coefficients is equal to 1).
Therefore, $u(P(5,-3,2)) \in [0, 2/3]$; but $s_3(P(5,-3,2)) = 1$.
Next, the linear independence of $\{\tau, s_2, s_3\}$ has been proven in (i). So it is
enough to show the existence of a knot $K$ with $\tau(K) = s_2(K) = s_3(K) = 0$, but $s_N(K) \neq 0$.
For this purpose, take
\[
K = P (5 , - 3 , 2) \# P (5 , - 3 , 2) \# - P (7 , - 5 , 4) \# - T (3 , 2).
\]
We have $s_2(K) = \tau(K) = 0$, and using FoamHo \cite{foamho}, $s_3(P(5,-3,2) = 1 \implies s_3(K) = 0$.
On the other hand, $s_N(P(5,-3,2)) \leq 2/(N-1) \implies s_N(K) \leq 4/(N-1) - 2 \neq 0$ since $N \geq 4$.
\end{Pf}
Let us compute another example, to illustrate that the Rasmussen invariant does not necessarily
give the best slice genus bound among the \Sln-concordance invariants.
\begin{Ex}
Let $K = 12n_{340}$, then $s_2(K) = 0, s_3(K) = 1$ and for $N \geq 4: s_N(K) \in \{2-2/(N-1),2\}$.
\end{Ex}
\begin{Pf}
The value of $s_2$ and $s_3$ may be computed using JavaKh \cite{javakh} and FoamHo \cite{foamho}, respectively;
the other values can be read from $\Hom[r]{K}{\infty}$, which we are going to compute using the
skein long exact sequence.
Notice that the calculation is rather quick, and that we do not need to determine the \homf-homology
of $K$ completely (this would be possible though,
using Rasmussen's spectral sequences \thmref{prop:Rasmussens-spectral-sequence}).

Resolving the crossing indicated in \thmref{fig:12n340} gives
$K$ as $K_+$, $10_{141}$ as $K_-$ and the positive Hopf link as $L_0$.
Resolving once more the indicated crossing of $10_{141}$ gives $10_{141}$ as $K_+$,
$8_9$ as $K_-$ and the positive Hopf link as $L_0$.
The knot $8_9$ is two-bridge, so its reduced \homf-homology is determined by its
\homf-polynomial and its signature. One finds
\[
\xdim\Hom[r]{8_9}{\infty}^{-4} = 
q^4 a^2.
\]
Applying \thmref{prop:homf les} twice gives
\[
\xdim\Hom[r]{K}{\infty} \leq t^4a^{-4} \xdim\Hom[r]{8_9}{\infty} + (t^{1/2}a^{-1} + t^{5/2}a^{-3})\xdim\Hom[rr]{T(2,2)}{\infty},
\]
and therefore
\[
\xdim\Hom[r]{K}{\infty}^0 \leq q^4a^{-2} + q^2a^{-2}.
\]
By \thmref{cor:tool} it follows that
\[
\forall N \geq 2: s_N(K) \in \{2-2/(N-1),2\}.
\]
\end{Pf}